\documentclass[12pt]{article}
\usepackage{bm}
\usepackage{amsmath}
\usepackage{amssymb}
\usepackage[mathscr]{eucal}
\usepackage{graphicx}
\usepackage{color}
\usepackage{multirow}
\usepackage{booktabs}

\newtheorem{Theorem}{\sc Theorem}

\newtheorem{Lemma}[Theorem]{\sc Lemma}
\newtheorem{Corollary}[Theorem]{\sc Corollary}
\newtheorem{Remark}[Theorem]{\sc Remark}

\newtheorem{Problem}[Theorem]{\sc Problem}

\newcommand{\ignore}[1]{}
\newcommand{\tabincell}[2]{\begin{tabular}{@{}#1@{}}#2\end{tabular}}  

\newcommand{\cS}{\mbox{{${\cal S}$}}}
\newcommand{\R}{{\if mm {\rm I}\mkern -3mu{\rm R}\else \leavevmode
\hbox{I}\kern -.17em\hbox{R} \fi}}

\newcounter{theorem}

\def\sqr#1#2{{
    \vcenter{
         \vbox{\hrule height.#2pt
               \hbox{\vrule width.#2pt height#1pt \kern#1pt
                     \vrule width.#2pt
               }
               \hrule height.#2pt
         }
    }
}}

\def\bar{\overline}

\def\lista#1
{{ \itemindent 0.0cm \labelsep .2cm \leftmargin 0.8cm \rightmargin
0.0cm \labelwidth 0.6cm \topsep 0.0mm
\parsep 0.0mm
\itemsep 0.0mm
\begin{list}{}
{ \setlength{\leftmargin}{.8cm} \setlength{\rightmargin}{0.0cm}
\setlength{\parsep}{0.0mm} \setlength{\topsep}{.0mm}
\setlength{\parskip}{.0cm} \setlength{\itemsep}{.0cm} }
{#1}\end{list}} }

\begin{document}

\title{\bf Numerical Analysis of History-dependent Variational-hemivariational Inequalities}


\author{
Shufen Wang\thanks{School of Mathematical Sciences, Fudan University, Shanghai 200433, China.},
\and Wei Xu\thanks{Tongji Zhejiang College, Jiaxing, Zhejiang 314051, China.}, \and
Weimin Han\thanks{Department of Mathematics, University of Iowa, Iowa City, IA
52242, USA.},\and Wenbin Chen\thanks{School of Mathematical Sciences and Shanghai Key Laboratory for
Contemporary Applied Mathematics, Fudan
University, Shanghai 200433, China. } }


\date{}
\maketitle


\noindent {\bf Abstract} \
In this paper, numerical analysis is carried out for a class of history-dependent variational-hemivariational inequalities arising in contact problems. Three different numerical treatments for temporal discretization are proposed to approximate the continuous model. Fixed-point iteration algorithms are employed to implement the implicit scheme and the convergence is proved with a convergence rate independent of the time step-size and mesh grid-size. A special temporal discretization is introduced for the history-dependent operator, leading to numerical schemes for which the unique solvability and error bounds for the temporally discrete systems can be proved without any restriction on the time step-size. As for spatial approximation, the finite element method is applied and an optimal order error estimate for the linear element solutions is provided under appropriate regularity assumptions. Numerical examples are presented to illustrate the theoretical results.

\vskip 2mm


\vskip 2mm










\section{Introduction}
\label{sec-1}
\setcounter{equation}0
The theory of variational and hemivariational inequalities plays an important role in the study of
nonlinear problems arising in Contact Mechanics, Physics, Economics and Engineering.  It is generally
agreed that interest in variational inequalities started with a contact problem posed by Signorini in 1930s.
The mathematical theory of variational inequalities relies on the properties of monotonicity, convexity and the subdiffierential of a convex function. Existence and uniqueness results can be found in
\cite{ LS1967, Bre1968, KS1980}. In terms of the numerical analysis for variational inequalities, the readers are referred to, e.g., \cite{GLT1981, Gl1984, HHNL1988}. Hemivariational inequalities as a useful generalization of variational inequalities were introduced in early 1980s by Panagiotopoulos (\cite{Pa83}). For hemivariational
inequalities, the notion of the subdifferential of in the sense of Clarke (\cite{Cla75, Cla1983}), defined for locally Lipschitz function, plays an important role.   Mathematical theory of hemivariational inequalities
is documented in several research monographs, e.g., \cite{Pa1993, NP1995, CLM2007, MOS2013, SM2018}.
A comprehensive reference on the numerical solution of hemivariational inequalities is \cite{HMP1999}
where the finite element method is applied to solve hemivariational inequalities, convergence of the
numerical solution is discussed, and solution algorithms are proposed and tested.  More recently, there
has been extensive research effort on optimal order error estimation and general convergence analysis of
numerical solutions for hemivariational inequalities, e.g., \cite{HMS14, BBHJ15, HSB17, HSD18, Han18},
and the survey paper \cite{HS19}.

Variational-hemivariational inequalities are a particular family of hemivariational inequalities,
having a special structure that include both convex and nonconvex functionals. Such inequalities arise
naturally in mathematical models for many contact problems, see \cite{SM2018} and the references therein.
A class of history-dependent variational-hemivariational inequalities with convex constraint
is studied in~\cite{SM16}. The novel structure of the inequalities involves a history-dependent
operator, unilateral constraint and two nondifferential functions, one of which is convex and the
other may be nonconvex. Existence, uniqueness and continuous dependence results are shown on the
inequalities, and are applied to the study of a quasistatic frictionless contact problem.
Numerical approximations of the history-dependent variational-hemivariational inequalities are
the topic of \cite{XHHCW19a}, where the second order accuracy for temporal discretization is achieved
by using the trapezoidal rule to approximate the history-dependent term. The spatial discretization is
done using the linear finite element and an optimal order error estimate is proved. Note that for the
numerical method studied in \cite{XHHCW19a}, a restriction on the time step-size is needed to
ensure the unique solvability of the numerical solution. In this paper, we develop new numerical methods to
solve the history-dependent variational-hemivariational inequalities with the property that no
restriction on the time step-size is needed for the unique solvability of the numerical solution.
Specifically, we use a partial trapezoidal rule to approximate the history-dependent operator, i.e., we
modify the trapezoidal rule by applying the left-point rectangular rule for the sub-integral over
the last time sub-interval. Consequently, the history dependent term is treated explicitly without loss of accuracy. This explicit treatment of the history dependent term eliminates the need for a restriction on
the time step-size. Although the explicit treatment is given in history dependent term, other implicit terms in the numerical scheme remain. We provide a fixed-point iterative algorithm to implement the implicit scheme and prove convergence of the iterative scheme, with a convergence rate independent of the time step-size and
the mesh grid-size. 
In addition, we propose two more schemes to solve the history-dependent variational-hemivariational
inequalities. One is of first order and the other is of second order with a slightly stringent small condition compared to that of the other two schemes. For all the three schemes, optimal order error estimates
with linear finite elements for spatial approximation are shown.

The rest of the paper is organized as follows. In Section \ref{sec-2},  we review some preliminary
material on functional analysis and present the history-dependent variational-hemivariational
inequality problem. In Section \ref{sec-3}, we propose three temporally semi-discrete schemes to approximate the continuous problem and error estimates are established. The corresponding fully discrete schemes are provided in Section \ref{sec-4}, and the error estimates are derived for the discrete problems with or without convex constraints. To implement the second order implicit scheme, in Section \ref{sec-6} we describe a fixed-point iterative process and prove that the iteration converges linearly with a convergence rate independent of the
time step-szie and mesh grid-size. Then in Section \ref{sec-5} we apply the theoretical results developed
in the previous sections in the numerical solution of a viscoelastic contact problem and obtain an optimal
order error estimate for the linear finite element
solutions under appropriate solution regularity assumptions. In Section \ref{sec:numer test} we report
results from simulation tests, focusing on the numerical evidence of the convergence orders.

\section{Preliminaries}   \label{sec-2}
\setcounter{equation}0

In this section we recall some notation, definitions and preliminary materials. Then we present
a class of history-dependent variational-hemivariational inequalities introduced in \cite{SM16}.

For normed spaces $X$ and $X_j$, let $X^*$ and $X_j^*$ be their topological duals, and write
$\|\cdot\|_X$, $\|\cdot\|_{X_j}$, $\|\cdot\|_{X^*}$ and $\|\cdot\|_{X_j^*}$ for their norms.
The duality pairing between $X$ and $X^*$, $\langle \cdot,\cdot\rangle_{X^* \times X}$,
is usually simply written as $\langle \cdot,\cdot\rangle$.  Similarly, the duality pairing
between $X_j^*$ and $X_j$, $\langle \cdot,\cdot\rangle_{X_j^* \times X_j}$,
is usually written as $\langle \cdot,\cdot\rangle_{X_j}$.

For a convex function $\varphi:X\rightarrow \mathbb{R}\cup \{+\infty\}$, the subset $\partial\varphi(x)$
of $X^*$,
\[ \partial\varphi(x)=\{x^*\in X^*\mid\varphi(v)-\varphi(x)\ge\langle x^*,v-x\rangle_{X^*\times X}\
\forall\,v\in X \} \]
is called the subdifferential (\cite{Rock1970}) of $\varphi$. If $\partial\varphi(x)$ is non-empty,
any element $x^* \in \partial\varphi(x)$ is called a subgradient of $\varphi$ at $x$ .
Let $\phi:X\rightarrow \mathbb{R}$ be a locally Lipschitz function. The generalized (Clarke) directional
derivative of $\phi$ at $x$ in the direction $v\in X$ is defined by (cf.\ \cite{Cla1983})
\[ \phi^{0}(x;v)=\limsup_{y\rightarrow x,\,\lambda \downarrow 0} \frac{\phi(y+\lambda v)-\phi(y)}{\lambda}. \]
The generalized gradient (subdifferential) of $\phi$ at $x$ is a subset of the
dual space $X^*$ given by
\[ \partial\phi(x)=\{\xi\in X^*\mid \phi^{0}(x;v)\ge\langle\xi,v\rangle_{X^*\times X}\ \forall\,v\in X \}. \]
An operator $A: X \rightarrow X^*$ is pseudomonotone (\cite{MOS2013}) if it is bounded and $u_n \rightarrow u$
weakly in $X$ together with $\limsup_n \langle Au_n, u_n-u \rangle_{X^* \times X} \leq 0$ imply
\[ \langle Au,u-v\rangle_{X^*\times X}\le\liminf_n\langle Au_n,u_n-v\rangle_{X^*\times X}\quad\forall\,v\in X. \]

Next we turn to some preliminary materials on function spaces and related operators. Following the standard
notation, we denote by $\mathbb{N}$ the set of positive integers, $\mathbb R_{+}=[0, +\infty)$ the
set of nonnegative real numbers, $C(\mathbb R_{+};X)$  and $C^{1}(\mathbb R_{+};X)$ the spaces
of continuous and continuously differentiable functions from $ \mathbb{R}_+ $ to $X$, respectively.
It is well known that if $X$ is a Banach space, $C(\mathbb R_{+};X)$ can be organized in a canonical
way as a Fr\'{e}chet space, i.e., it is a complete metric space in which the corresponding topology
is induced by a countable family of seminorms. Furthermore,
$x_k\rightarrow x$ in $C(\mathbb{R}_{+}; X)$ as $k\rightarrow \infty$ if and only if
$\max\limits_{r\in [0,n]}\|x_k(r)-x(r)\|_X\rightarrow 0$ as $k\rightarrow \infty$ for all $n\in\mathbb{N}$.

Let there be given two normed spaces $X$ and $Y$.  Following \cite{SM11}, an operator
$\cS:C(\mathbb R_{+};X)\rightarrow C(\mathbb R_{+};Y)$ is called history-dependent if for any
$n\in \mathbb N$, there exists an $s_n >0$ such that for all $t\in[0,n]$,
\begin{equation}
\|(\cS u_1)(t)-(\cS u_2)(t)\|_Y \leq s_n \int_0^t \|u_{1}(s)-u_{2}(s)\|_{X}ds
\quad \forall\, u_{1}, u_{2} \in C(\mathbb R_{+};X).
\label{pre-1}
\end{equation}

Now we are in a position to introduce the variational-hemivariational inequalities.
Let $X$, $X_j$, $Y$ be normed spaces and $K\subset X$. Given operators $A:X \rightarrow X^*$,
$\cS:C(\mathbb{R}_+;X) \rightarrow C(\mathbb{R}_+;Y)$, $\gamma_j: X \rightarrow X_j$ and functions
$\varphi:Y\times K \times K \rightarrow \mathbb{R}$, $j:X_j \rightarrow \mathbb{R}$, we consider the following problem (\cite{SM16,XHHCW19a}).

\begin{Problem} \label{Problem 1}
Find $u \in C(\mathbb R_{+};K)$ such that for all $t\in\mathbb{R_+}$,
\begin{equation}
    \begin{aligned}
  & \langle Au(t),v-u(t) \rangle + \varphi((\cS u)(t),u(t),v)- \varphi((\cS u)(t),u(t),u(t))\\
    & \quad{}+j^0(\gamma_{j}u(t);\gamma_{j}v-\gamma_{j}u(t))\ge\langle f(t),v-u(t)\rangle\quad\forall\,v\in K.
    \end{aligned}
    \label{pre-2}
\end{equation}
\end{Problem}

In the study of Problem \ref{Problem 1}, the following hypotheses are adopted (\cite{SM16,XHHCW19a}):
\begin{eqnarray}
&& \begin{array}{ll}
	\ X \textrm{ is a reflexive Banach space},\ K \textrm{ is a closed and convex subset of} \ X \\
	\textrm{with} \ 0\in K.
	\end{array}
   \label{pre-3}\\[1mm]
&&  \left \{
     \begin{array}{ll}
      X_{j}\ \textrm{is a Banach space}, \ \gamma_j \in \mathcal{L}(X;X_j), \ \textrm{there exists}\ c_{j}>0 \ \textrm{such that} \\
      \qquad \qquad \|\gamma_{j}v\|_{X_j} \leq c_{j}\|v\|_X \quad \forall \, v \in X.\\
     \end{array}
   \right.
   \label{pre-4}\\ [1mm]
&& \left \{
     \begin{array}{ll}
       A: X \rightarrow X^* \ \textrm{is  an  operator  such that} \\
       \quad (a)\ A\ \textrm{is Lipschitz continuous with a Lipschitz constant}\ L_A>0.\\
       \quad (b) \ A \ \textrm{is  strongly monotone, i.e., there exists} \  m_A > 0 \ \textrm{such that}\\
       \qquad \langle Av_1-Av_2,v_1-v_2\rangle \geq m_A\|v_1-v_2\|_{X}^{2} \quad \forall\, v_1, v_2 \in X.\\
     \end{array}
   \right.
   \label{pre-5}\\[1mm]
&& \left \{
     \begin{array}{ll}
       \varphi:Y\times K\times K\rightarrow \mathbb{R} \ \textrm{is a function such that}\\
       \quad \textrm{(a)}\  \varphi(y,u,\cdot):K\rightarrow \mathbb{R} \ \textrm{is convex and l.s.c on } K,  \forall\, y \in Y, \forall\, u \in K.\\
       \quad \textrm{(b)} \ \textrm{there exists} \  \alpha_{\varphi} >0 \ \textrm{and} \  \beta_{\varphi}> 0 \ \textrm{such that}\\
       \qquad  \varphi(y_1,u_1,v_2)-\varphi(y_1,u_1,v_1)+\varphi(y_2,u_2,v_1)-\varphi(y_2,u_2,v_2)\\
       \qquad\quad\le\alpha_{\varphi}\|u_1-u_2\|_X\|v_1-v_2\|_X+\beta_{\varphi}\|y_1-y_2\|_{Y}\|v_1-v_2\|_X\\
       \qquad  \quad   \forall\, y_1, y_2\in Y,\ \forall\, u_1, u_2, v_1, v_2 \in K.
     \end{array}
   \right.
   \label{pre-7} \\[1mm]
   &&
   \begin{array}{ll}
   \cS:C(\mathbb R_{+};X)\rightarrow C(\mathbb R_{+};Y)\  \textrm{is a history-dependent operator}.
   \end{array}
   \label{pre-6}  \\[1mm]
&& \left \{
     \begin{array}{ll}
     j: X_{j} \rightarrow  \mathbb R\   \textrm{is a function such that}\\
     \quad \textrm{(a)} \ j \ \textrm{is locally Lipschitz}.\\
     \quad \textrm{(b)} \  \|\partial j(z)\|_{X^*_j}\leq c_0+c_1 \|z\|_{X_j}\quad \forall\, z \in X_j \ \textrm{with}\ c_0, c_1\geq 0.\\
     \quad \textrm{(c)} \ \textrm{there exists}\ \alpha_j >0 \ \textrm{such that}\\
     \qquad j^{0}(z_1;z_2-z_1)+j^{0}(z_2;z_1-z_2)\le\alpha_j\|z_1-z_2\|_{X_j}^2\quad\forall\,z_1, z_2\in X_j.\\
     \end{array}
   \right.
   \label{pre-8} \\[1mm]
&& \qquad \qquad \qquad \qquad \qquad f \in C(\mathbb R_{+};X^*). \label{pre-9} \\[1mm]
&& \qquad \qquad \qquad \qquad \qquad \alpha_\varphi + \alpha_j c_j^2 < m_A. \label{pre-10}
\end{eqnarray}
The space $X_{j}$ is introduced for convenience of error estimation for the discrete problems.
For a specific contact problem, $X_j$ can be the space of square integrable functions over the contact boundary
and $\gamma_j: X \rightarrow X_j$ is the corresponding trace operator. For a locally Lipschitz function $j$,
\eqref{pre-8}\,(c) is equivalent to the following relaxed monotonicity condition
\[\langle \partial j(z_1)-\partial j(z_2), z_1-z_2 \rangle \geq -\alpha_j \|z_1-z_2\|_{X_j}^2
\quad\forall\, z_1, z_2 \in {X_j}. \]

The unique solvability of Problem \ref{Problem 1} has been shown in \cite{XHHCW19a} under the conditions
\eqref{pre-3}--\eqref{pre-10}. We will consider the following form of
the operator $\cS: C(\mathbb R_{+};X) \rightarrow C(\mathbb R_{+};Y)$ (\cite{SM11}):
\begin{equation}
(\cS v)(t)=R\left(\int_{0}^{t}q(t,s)v(s)ds+a_{S}\right)\quad\forall\,v\in C(\mathbb R_{+};X),\ \forall\,t\in\mathbb R_{+},
\label{pre-11}
\end{equation}
where $R \in \mathcal{L}(X;Y)$, $q \in C(\mathbb R_{+}\times \mathbb R_{+};\mathcal{L}(X))$, $a_S \in X$.
It can be shown that the operator $\cS$ given by (\ref{pre-11}) is a history-dependent operator.

\section{Temporally Semi-Discrete Approximations} \label{sec-3}
\setcounter{equation}0

In \cite{XHHCW19a}, a second-order numerical scheme is provided to approximate the continuous Problem \ref{Problem 1} with a restriction on the time step-size. In this section, we handle the history-dependent term in a different manner, and propose three temporally discrete schemes for solving Problem \ref{Problem 1} without any restriction on the time step-size.  Moreover, we derive the corresponding convergence results. Below we use $C$ to represent a positive constant independent of time step-size and mesh grid-size.
We use the standard notation for Sobolev spaces (cf.\ \cite{Adam2003}).

For a fixed $T\in \mathbb{R}_+$, we split the time interval $I=[0,T]$ by uniform partitions. Given a positive
integer $N$, let $k=T/N$ be the time step-size, and denote by $t_n = nk$, $0\le n\le N$, the nodes.  We
comment that all the discussions below can be extended to the case with non-uniform partitions of the
time interval.  For a  continuous function $v$ of the temporal variable $t$, we write $v_j=v(t_j)$,
$0\le j\le N$. For a discretization of the history dependent operator $\cS$ in (\ref{pre-11}), we employ
a modified trapezoidal rule to approximate the integral $\int_{0}^{t_n}q(t,s)v(s)ds $ in the sense that on
the last sub-interval $[t_{n-1},t_n]$, the left-point rectangular rule is applied. Recall the trapezoidal rule
\begin{equation}
    \int_{0}^{t_n}Z(s)ds\approx  \frac{k}{2}Z(t_0)+k\sum_{j=1}^{n-1}Z(t_j)+\frac{k}{2}Z(t_n).
  \label{tem-1}
\end{equation}
The approximation of $\cS_n:=\cS(t_n)$ can be defined as follows:
\begin{equation}
\cS_{n,L}^{k}v:=R\left(\frac{k}{2}q(t_n,t_0)v_0+k\sum_{j=1}^{n-1}q(t_n,t_j)v_j+\frac{k}{2}q(t_n,t_{n-1})v_{n-1}+a_S\right).
\label{tem-2}
\end{equation}
Using arguments similar to that in \cite[Section 3]{KBHS14},  for $v\in W^{1,\infty}_{loc}(\mathbb{R_{+}};X)$
and $q\in C^1(\mathbb{R_+} \times \mathbb{R_+};\mathcal{L}(X))$, we have
\begin{equation} \label{eq:S accurancy1}
\|\cS_{n,L}^k v - \cS_n v\| \le Ck \|v\|_{W^{1,\infty}(I;X)},
\end{equation}
and for $v\in W^{2,\infty}_{loc}(\mathbb{R_{+}};X)$ and $q\in C^2(\mathbb{R_+}\times\mathbb{R_+};\mathcal{L}(X))$,
\begin{equation} \label{eq:S accurancy}
	\|\cS_{n,L}^k v - \cS_n v\| \le C\,k^2 \|v\|_{W^{2,\infty}(I;X)}.
\end{equation}

\begin{Remark} \label{2nd-operator def}
The choice of the operator $\cS_{n,L}^k$ used to approximate $\cS_{n}$ is not unique. For example, we may choose
\[ \hat{\cS}_{n}^kv :=R\left(\frac{k}{2}q(t_n,t_0)v_0+k\sum_{j=1}^{n-1}q(t_n,t_j)v_j
+\frac{k}{2}(2q(t_n,t_{n-1})v_{n-1}-q(t_n,t_{n-2})v_{n-2})+a_S\right) \]
which defines another second-order accurate approximation of $\cS_{n}$, or choose
\[ \widetilde{\cS}_{n}^kv :=R\left(k\sum_{j=0}^{n-1}q(t_n,t_j)v_j+a_S\right) \]
which is a first-order accurate approximation.
\end{Remark}
We note that the following weak formulation is equivalent to Problem \ref{Problem 1}:
\begin{Problem} \label{Problem equiv}
	Find $u \in C(\mathbb{R_{+}}; K)$ such that for all $t\in\mathbb{R_+}$
	\begin{equation} \label{eq: equiv}
	\begin{aligned}
		  & \langle Au(t),v-u(t) \rangle + \varphi((\cS u)(t),u(t),v)- \varphi((\cS u)(t),u(t),u(t))  \\
		& \qquad{}+j^0(\gamma_{j}u(t);\gamma_{j}v-\gamma_{j}u(t)) + \langle j_c(\gamma_ju(t)), \gamma_jv-\gamma_j u(t) \rangle_{X_j}   \\
		& \quad \ge\langle f(t),v-u(t)\rangle +  \langle j_c(\gamma_ju(t)), \gamma_jv-\gamma_j u(t) \rangle_{X_j} \quad\forall\,v\in K.
	\end{aligned}
	\end{equation}
\end{Problem}
In \cite{XHHCW19a}, $j_c$ is chosen as the differential of a quadratic function $\frac{\alpha}{2}\|u\|^2_{X_j}$. In this paper, we discuss about $j_c$ in a more general framework. Assume
\begin{equation}
	\left \{
	\begin{aligned}
		&j_c : X_j \rightarrow X_j^* \  \textrm{ is a linear operator such that} \\
		&\quad\textrm{(a)} \ \|j_c(z)\|_{X_j^*}\leq \alpha_{c}\|z\|_{X_j}\quad \forall\, z\in X_j; \\
		&\quad\textrm{(b)} \ \langle j_c(z),z\rangle_{X_j}\ge\alpha_{j}\|z\|^2_{X_j} \quad \forall \, z\in X_j.\\
	\end{aligned}
	\right.
	\label{pre-jc}
\end{equation}
The operator $j_c$ can be regarded as a convexification of $j^0$ in the sense that
\begin{equation} \label{eq:j_convex}
\begin{aligned}
& j^0(z_1;z_2-z_1)+j^0(z_2;z_1-z_2)+\langle j_c(z_1),z_2-z_1\rangle_{X_j}+\langle j_c(z_2),z_1-z_2\rangle_{X_j}\\
& \qquad \le \alpha_j\|z_1-z_2\|_{X_j}^2 +  \langle j_c(z_1-z_2), z_2-z_1 \rangle_{X_j} \\
& \qquad \le 0,
\end{aligned}
\end{equation}
where the last equality follows from \eqref{pre-jc} (b).


\subsection{A first-order temporally semi-discrete scheme}

The first order temporally semi-discrete scheme for Problem \ref{Problem 1} is the following.

\begin{Problem} \label{Problem 4}
Find a discrete solution $u^{k}:=\{u_{n}^{k}\}_{n=0}^{N}\subset K$ such that
\begin{equation}
  \begin{aligned}
  & \langle Au_{n}^{k},v-u_{n}^{k} \rangle+\varphi(\cS_{n,L}^{k}u^{k},u_{n-1}^{k},v)-\varphi(\cS_{n,L}^{k}u^{k},u_{n-1}^{k},u_{n}^{k}) \\
  & \qquad+j^0(\gamma_{j}u_{n}^{k};\gamma_{j}v-\gamma_{j}u_{n}^{k}) +\langle j_c(\gamma_ju_n^k),\gamma_jv-\gamma_ju_n^k\rangle_{X_j} \\
  &\quad \geq \langle f_n,v-u_{n}^{k}\rangle{ +\langle j_c(\gamma_ju_{n-1}^k),\gamma_jv-\gamma_ju_n^k\rangle_{X_j}} \quad \forall\, v \in K.
  \end{aligned}
  \label{tem-3}
\end{equation}
\end{Problem}

\begin{Remark} \label{remark:1rdexist}
Note that the approximation $\cS_{n,L}^ku^k$ for the history dependent operator does not involve information
on the current numerical solution $u_n^k$, and the second argument of $\varphi$ is explicitly treated, which is important for  numerical implementation.
The function $\varphi$ appeared in \eqref{tem-3} is convex with respect to the unknown variable (the third
argument) according to assumption \eqref{pre-7}. Moreover, $j_c$ plays the role to convexify the function $j$,
i.e., $j^0(\gamma_{j}u_{n}^{k};\gamma_{j}v-\gamma_{j}u_{n}^{k}) +\langle j_c(\gamma_ju_n^k),\gamma_jv-\gamma_ju_n^k\rangle_{X_j}$ becomes the directional derivative of a convex function. Therefore, convex optimization techniques could be applied to solve the inequality \eqref{tem-3} and the unique solvability of Problem
\ref{Problem 4} can be obtained without the constraint \eqref{pre-10} by applying results on
elliptic variational-hemivariational inequality (\cite{MOS17}). Specifically,  the operator $T_1$ defined by $T_1v=Av+\partial j(v)+j_c(v)$ is bounded, coercive and pseudomonotone, the function $\varphi(v)$ can be extended to $X$, denoted as $\widetilde{\varphi}(v)$ with $\widetilde{\varphi}(v)=+\infty$ for any $v\in X\backslash K$. In this way, the operator $T_2$ with $T_2 v = \partial \widetilde{\varphi}(v)$ is maximal monotone.
Hence, Problem \ref{Problem 4} has a unique solution.
\end{Remark}

\begin{Remark} \label{jc_choice}
The choice of $j_c$ is not unique.  The critical point is that $j_c$ should be ``convex'' enough to have the
non-convexity of $j^0$ under control, i.e., the inequality \eqref{eq:j_convex} is required.  On the other
hand, we can split $j^0$ in another way, e.g.,
\[ \langle j_c(\gamma_ju_n^k),\gamma_jv-\gamma_ju_n^k\rangle_{X_j}
+ \big(j^0(\gamma_{j}u_{n-1}^{k};\gamma_{j}v-\gamma_{j}u_{n}^{k})
- \langle j_c(\gamma_ju_{n-1}^k),\gamma_jv-\gamma_ju_n^k\rangle_{X_j}\big)\]
could be used to approximate $j^0(\gamma_{j}u_{n}^{k};\gamma_{j}v-\gamma_{j}u_{n}^{k})$. In this way,
the inequality \eqref{tem-3} becomes a convex problem with linear operators, for which efficient numerical algorithms are available.
\end{Remark}

According to the statement in Remark \ref{remark:1rdexist}, we have the following unique solvability result for Problem \ref{Problem 4}.

\begin{Theorem} \label{thm:1rdunique}
Under the conditions~\eqref{pre-3}--\eqref{pre-9} and \eqref{pre-jc}, the semi-discrete Problem~\ref{Problem 4} is uniquely solvable.
\end{Theorem}

For error estimation, we first introduce some auxiliary techniques.

\begin{Lemma}\label{Lemma Gronwall}
Let $\{a_n\}$ be a nonnegative sequence satisfying
\begin{align*}
		a_{n} \le b_0+C_1k\sum_{j=0}^{n-1}a_j + \theta_1a_{n-1}+\theta_2a_{n-2} \quad \forall\, n\ge 2,
\end{align*}
where $a_0,a_1$, $b_0,\theta_1,\theta_2,C_1$ are nonnegative constants and $0\le\theta_1+\theta_2<1$.  Then
\begin{equation} \label{eq: grown1}
a_n\le \bigg(\frac{b_0}{1-\theta_1-\theta_2}+\frac{C_1k(a_0+a_1)}{1-\theta_1-\theta_2}+\theta_1a_1+\theta_2a_0\bigg)\bigg(1+\frac{C_1k}{1-\theta_1-\theta_2}\bigg)^{n-2}.
\end{equation}
\end{Lemma}
{\it Proof.} For convenience, let
\[\bar\alpha :=\frac{b_0}{1-\theta_1-\theta_2}+\frac{C_1k(a_0+a_1)}{1-\theta_1-\theta_2}+\theta_1a_1+\theta_2a_0.\]
We prove the result with an induction. 	For $n=2$, we have the following bound:
	\begin{align*}
a_{2} \le b_0+C_1k(a_1+a_0) +\theta_1a_1+\theta_2a_0\le \bar\alpha.
	\end{align*}
Thus, \eqref{eq: grown1} holds for $n=2$. Assume that for $n\le m$,
	 \begin{align*}
	 a_{n}\le \bar\alpha\big(1+\frac{C_1k}{1-\theta_1-\theta_2}\big)^{n-2}.
	 \end{align*}
Then for $n=m+1$,
\begin{align*}
a_{m+1}&\le b_0+C_1k\sum_{j=0}^{m}a_j+\theta_1a_m+\theta_2a_{m-1} \\
&\le b_0+C_1k\big[a_0+a_1+\bar\alpha\sum_{j=2}^{m}(1+\frac{C_1k}{1-\theta_1-\theta_2})^{j-2}\big]\\
&\quad{}+(\theta_1+\theta_2)\bar\alpha(1+\frac{C_1k}{1-\theta_1-\theta_2})^{m-2}\\
&=b_0+C_1k(a_0+a_1)+\bar\alpha\cdot C_1k\frac{(1+\frac{C_1k}{1-\theta_1-\theta_2})^{m-1}-1}{C_1k}(1-\theta_1-\theta_2)\\
&\quad{}+(\theta_1+\theta_2)\bar\alpha(1+\frac{C_1k}{1-\theta_1-\theta_2})^{m-2}\\
	 &\le \bar\alpha(1+\frac{C_1k}{1-\theta_1-\theta_2})^{m-1}(1-\theta_1-\theta_2)+(\theta_1+\theta_2)\bar\alpha(1+\frac{C_1k}{1-\theta_1-\theta_2})^{m-2}\\
&\le\bar\alpha(1+\frac{C_1k}{1-\theta_1-\theta_2})^{m-2}(1+C_1k)\le\bar\alpha(1+\frac{C_1k}{1-\theta_1-\theta_2})^{m-1},
	 \end{align*}
	 where we use the fact that
	 $$b_0+C_1k(a_0+a_1)-(1-\theta_1-\theta_2)\bar\alpha\le 0.$$	
This completes the proof.
	\hfill$\Box$
	
\begin{Corollary} \label{Col Gronwall}
		Assume that $\{a_n\}$ is a nonnegative sequence satisfying
		\begin{align*}
		a_{n} \le b_0+C_2k\sum_{j=0}^{n-1}a_j + \theta_1a_{n-1} \quad \forall\, n\ge 1,
		\end{align*}
where $a_0$, $b_0$, $\theta_1$ and $C_2$ are nonnegative constants and $\theta_1<1$.  Then
\begin{equation}
a_n\le \bigg(\frac{b_0}{1-\theta_1}+\frac{C_2k}{1-\theta_1}a_0+\theta_1a_0\bigg)\bigg(1+\frac{C_2k}{1-\theta_1}\bigg)^{n-1}.
\end{equation}
\end{Corollary}
\begin{Lemma} \label{Lem: square}
	Assume $e_0, e_1$ and $e_2$ are nonnegative numbers such that
	\begin{equation} \label{eq:square cond}
		 e_0^2 \le e_1e_0 + e_2^2,
	\end{equation}
	then
	\begin{equation} \label{eq: square}
		e_0 \le e_1 + e_2.
	\end{equation}
\end{Lemma}
{\it Proof.}
From \eqref{eq:square cond}, we have
\begin{equation}
	\left(e_0 - \frac{e_1}{2}\right)^2 \le \frac{e_1^2}{4} + e_2^2 \le \left(\frac{e_1}{2}+e_2\right)^2.
\end{equation}
Taking the square root of both sides gives \eqref{eq: square}.
	\hfill$\Box$

We now turn to an error analysis for Problem \ref{Problem 4}. For convenience, we denote
$\|R\|=\|R\|_{\mathcal{L}(X;Y)}$ and $\|q\|=\|q\|_{C(I\times I;\mathcal{L}(X))}$. The following
smallness condition is needed instead of the original one \eqref{pre-10}:
\begin{align} \label{eq: small condition 1rd}
\alpha_{\varphi} + \alpha_{c}c_j^2 < m_A.
\end{align}
\begin{Theorem} \label{Theorem 6}
	Assume \eqref{pre-3}--\eqref{pre-9}, \eqref{pre-jc}, \eqref{eq: small condition 1rd} and the regularity
		$q\in C^1(\mathbb R_{+} \times \mathbb R_{+};\mathcal{L}(X))$, $u\in {W_{loc}^{1,\infty}(\mathbb R_{+};X)}$.
Then for the semi-discrete solution of Problem \ref{Problem 4}, the following error bound holds:
	\begin{equation}
	\max \limits_{n\leq N}\|u_n-u_n^k\|_{X} \leq C_3k,
	\label{tem-oneorder-1}
	\end{equation}
	where $C_3>0$ is a constant independent of $k$.
\end{Theorem}
{\it Proof.}
We take $t=t_n$ in the inequality~\eqref{pre-2} to get
\begin{equation}
\begin{aligned}
& \langle Au_n,v-u_n \rangle + \varphi(\cS_{n}u,u_n,v)- \varphi(\cS_{n}u,u_n,u_n)\\
& \quad +j^0(\gamma_{j}u_n;\gamma_{j}v-\gamma_{j}u_n) \geq \langle f_n,v-u_n \rangle \quad \forall\, v \in K,
\end{aligned}
\label{tem-oneorder-2}
\end{equation}
where $\cS_{n}u=R(\int_0^{t_n}q(t_n,s)u(s)ds+a_{S})$. Let $v= u_n^k$ in~\eqref{tem-oneorder-2},
\begin{equation}
\begin{aligned}
& \langle Au_n,u_n^k-u_n \rangle+\varphi(\cS_{n}u,u_n,u_n^k)-\varphi(\cS_{n}u,u_n,u_n)\\
& \quad +j^0(\gamma_{j}u_n;\gamma_{j}u_n^k-\gamma_{j}u_n) \geq \langle f_n,u_n^k - u_n \rangle.
\end{aligned}
\label{tem-oneorder-3}
\end{equation}
Taking $v=u_n$ in~\eqref{tem-3} yields
\begin{equation}
\begin{aligned}
& \langle A u_n^k, u_n-u_n^k \rangle + \varphi(\cS_{n,L}^k u^k, u_{n-1}^k, u_n)- \varphi(\cS_{n,L}^k u^k, u_{n-1}^k, u_n^k)\\
& \qquad\quad {}+j^0(\gamma_{j}u_n^k;\gamma_{j}u_n-\gamma_{j}u_n^k)+\langle j_c(\gamma_ju_n^k),\gamma_ju_n-\gamma_ju_n^k\rangle_{X_j} \\
&\qquad \ge\langle f_n,u_n-u_n^k \rangle+\langle j_c(\gamma_ju_{n-1}^k),\gamma_ju_n-\gamma_ju_n^k\rangle_{X_j}.
\end{aligned}
\label{tem-oneorder-4}
\end{equation}
Adding~\eqref{tem-oneorder-3} to~\eqref{tem-oneorder-4} and employing the strong monotonicity of $A$, we obtain
\begin{align*}
m_A\|u_n-u_n^k\|^2_X &\le \langle A u_n-A u_n^k, u_n-u_n^k \rangle \\
& \leq \varphi(\cS_{n}u,u_n,u_n^k)-\varphi(\cS_{n}u,u_{n},u_n)+\varphi(\cS_{n,L}^k u^k,u_{n-1}^k,u_n)\\
&  \quad{} -\varphi(\cS_{n,L}^k u^k,u_{n-1}^k,u_n^k)+j^0(\gamma_{j}u_n;\gamma_{j}u_n^k-\gamma_{j}u_n)\\
&  \quad{} +j^0(\gamma_{j}u_n^k;\gamma_{j}u_n-\gamma_{j}u_n^k)+\langle j_c(\gamma_ju_n^k),\gamma_ju_n-\gamma_ju_n^k\rangle_{X_j}\\
& \quad{} -\langle j_c(\gamma_ju_{n-1}^k),\gamma_ju_n-\gamma_ju_n^k\rangle_{X_j},
\end{align*}
which is rewritten as
\begin{equation}
m_A\|u_n-u_n^k\|^2_X\le E_{\varphi} + E_j + E_{j_c},
\label{tem-oneorder-5}
\end{equation}
where
	\begin{eqnarray}
	&& \begin{array}{ll}
	E_{\varphi}=\varphi(\cS_{n,L}^{k}u^{k},u_{n-1}^{k},u_n)-\varphi(\cS_{n,L}^{k}u^{k},u_{n-1}^{k},u_{n}^{k}) \\
	\qquad \quad+
	\varphi(\cS_{n}u,u_n,u_{n}^{k})-\varphi(\cS_{n}u,u_n,u_n),
	\end{array}
	\label{tem-Ephi-express}\\[1mm]
	&& \begin{array}{ll}
	E_{j_c}= \langle j_c(\gamma_ju_n),\gamma_ju_n-\gamma_ju_n^k\rangle_{X_j}-\langle j_c(\gamma_ju_{n-1}^k),\gamma_ju_n-\gamma_ju_n^k\rangle_{X_j},
	\end{array}
	\label{tem-Ejc-express} \\[1mm]
	&& \begin{array}{ll}
	E_j = j^0(\gamma_{j}u_{n};\gamma_{j}u_{n}^{k}-\gamma_{j}u_{n})+j^0(\gamma_{j}u_{n}^{k};\gamma_{j}u_n-\gamma_{j}u_{n}^{k})\\
	\qquad \quad +\langle j_c(\gamma_ju_n),\gamma_ju_n^{k}-\gamma_ju_n\rangle_{X_j}+\langle j_c(\gamma_ju_n^{k}),\gamma_ju_n-\gamma_ju_n^{k}\rangle_{X_j}.
	\end{array}
	\label{tem-Ej-express}
	\end{eqnarray}
The term $E_j$ can be bounded by zero from above according to \eqref{eq:j_convex}. Utilizing the regularity of $u$ and the properties of $j_c$ gives
	\begin{equation}
	\begin{aligned}
		E_{j_c}&=\langle j_c(\gamma_ju_n-\gamma_ju_{n-1}^k),\gamma_ju_n-\gamma_ju_n^k\rangle_{X_j} \\
		&=\langle j_c(\gamma_ju_n-\gamma_ju_{n-1}),\gamma_ju_n-\gamma_ju_n^k\rangle_{X_j} \\
		& \quad +\langle j_c(\gamma_ju_{n-1}-\gamma_ju_{n-1}^k),\gamma_ju_n-\gamma_ju_n^k\rangle_{X_j}\\
		&\le \alpha_c c_j^2\big(k\|u\|_{W^{1,\infty}(I,X)}+\|u_{n-1}-u_{n-1}^k\|_X\big)\|u_n-u_n^k\|_X.
	\end{aligned}
	\label{tem-Ejc}
	\end{equation}
From \eqref{pre-7} we can see
	\begin{equation}
	\begin{aligned}
		E_{\varphi} &\le \big(\alpha_\varphi \|u_n-u_{n-1}^k\|_X+\beta_\varphi\|\cS_{n}u-\cS_{n,L}^k u^k\|_Y\big) \|u_n-u_n^k\|_X \\
		&\le \big(k\alpha_\varphi\|u\|_{W^{1,\infty}(I,X)}+\alpha_\varphi \|u_{n-1}-u_{n-1}^k\|_X \\
		& \quad +\beta_\varphi\|\cS_{n}u-\cS_{n,L}^k u^k\|_Y\big) \|u_n-u_n^k\|_X.
	\end{aligned}
		\label{tem-Ephi}
	\end{equation}
From \eqref{eq:S accurancy1}, it holds
	\begin{equation}
	\begin{aligned}
	\|\cS_n u-\cS_{n,L}^k u^k\|_{Y} & \leq \|\cS_n u-\cS_{n,L}^k u\|_Y+ \|\cS_{n,L}^k u-\cS_{n,L}^k u^k\|_Y\\
	& \leq Ck\|u\|_{W^{1,\infty}(I;X)}+ \frac{3}{2}k\|R\|\|q\|\sum_{j=0}^{n-1}\|u_j-u_j^k\|_X.
	\end{aligned}
	\label{tem-oneorder-6}
	\end{equation}
From \eqref{tem-oneorder-5} and \eqref{tem-Ejc}--\eqref{tem-oneorder-6}, we obtain
	\begin{equation}
	\begin{aligned}
	m_A\|u_n-u_n^k\|_X\ &\le  Ck\|u\|_{W^{1,\infty}(I;X)} +\frac{3}{2}k\beta_{\varphi}\|R\|\|q\|\sum_{j=0}^{n-1}\|u_j-u_j^k\|_X \\
	& \quad +(\alpha_{\varphi}+\alpha_cc_j^2)\|u_{n-1}-u_{n-1}^k\|_X.
	\end{aligned}
	\label{tem-oneorder-7}
	\end{equation}
By applying Corollary \ref{Col Gronwall},
\begin{equation}
	\begin{aligned}
	\|u_n-u_n^k\|_X&\le \left(k\frac{C\|u\|_{W^{1,\infty}(I;X)}}{m_A-\alpha_\varphi-\alpha_cc_j^2}+\left(\frac{\frac32k\beta_\varphi\|R\|\|q\|}{m_A-\alpha_\varphi-\alpha_{c}c_j^2} \right.\right.\\
	& \left.\left. \qquad+\frac{\alpha_{\varphi}+\alpha_{c}c_j^2}{m_A}\right)\|u_0-u_0^k\|_X\right)\cdot \left(1+k\frac{\frac32\beta_\varphi\|R\|\|q\|}{m_A-\alpha_{\varphi}-\alpha_{c}c_j^2}\right)^{n-1}.
	\end{aligned}
\end{equation}
\noindent
Note that when $t=t_0=0$, the integral of history-dependent operator is zero and there is no temporally discrete error; thus $\|u_0-u_0^k\|_X=0$. Then,
\begin{align*}
	\|u_n-u_n^k\|_X&\le k\frac{C\|u\|_{W^{1,\infty}(I;X)}}{m_A-\alpha_\varphi-\alpha_cc_j^2}\cdot \bigg(1+k\frac{\frac32\beta_\varphi\|R\|\|q\|}{m_A-\alpha_{\varphi}-\alpha_{c}c_j^2}\bigg)^{n-1} \\
	& \le C_3 k,
\end{align*}
where
$$C_3 =  \frac{C\|u\|_{W^{1,\infty}(I;X)}}{m_A-\alpha_\varphi-\alpha_cc_j^2}\cdot \exp \bigg\{\frac{\frac32\beta_\varphi\|R\|\|q\|}{m_A-\alpha_{\varphi}-\alpha_{c}c_j^2}t_n\bigg\}, $$
and the error bound~\eqref{tem-oneorder-1} follows.
\hfill$\Box$

\begin{Remark} \label{remark:operatorS-1rd}
The first-order accuracy remains valid if $\widetilde{\cS}_n^k$ is used to approximate the history-dependent
operator $\cS$ in the temporally semi-discrete scheme \eqref{tem-3}.
\end{Remark}
\subsection{Second-order temporally semi-discrete schemes}

In this subsection, we propose and study two second-order schemes to temporally approximate
Problem \ref{Problem 1}.  The first scheme is the following.

\begin{Problem}\label{problem tem-impli-expli}
Find $u^k:=\{u_n^k\}_{n=0}^{N}\subset K$ such that
\begin{equation}
\begin{aligned}
& \langle Au_n^k,v-u_n^k\rangle+\varphi(\cS_{n,L}^k u^k,u_n^k,v)-\varphi(\cS_{n,L}^k u^k,u_n^k,u_{n}^{k})\\
& \qquad\quad{}+j^0(\gamma_{j}u_n^k;\gamma_j v-\gamma_j u_n^k)\ge\langle f_n,v-u_n^k\rangle\quad\forall\,v \in K.
\end{aligned}
\label{tem-impli-expli}
\end{equation}
\end{Problem}

Note that the history-dependent operator is approximated using available numerical solution values and
the current unknown value $u_n^k$ is not involved. In this way, unlike the numerical scheme studied
in \cite{XHHCW19a}, the semi-discrete Problem~\ref{problem tem-impli-expli} is ensured to have
a unique solution regardless of the size of the time step-size using the same Banach fixed-point argument
as in \cite{XHHCW19a}.
	
\begin{Theorem} \label{thm:2rdunique}
Under the conditions~\eqref{pre-3}--\eqref{pre-10}, the semi-discrete Problem~\ref{problem tem-impli-expli} has a unique solution.
\end{Theorem}

We turn to the error estimation of Problem \ref{problem tem-impli-expli}.

\begin{Theorem}\label{Theorem tem-impli-expli}
		Assume \eqref{pre-3}--\eqref{pre-10} and the regularity
		$q\in C^2(\mathbb R_{+} \times \mathbb R_{+};\mathcal{L}(X))$, $u\in {W_{loc}^{2,\infty}(\mathbb R_{+};X)}$. Then for the semi-discrete solution of Problem \ref{problem tem-impli-expli}, we have the error bound
\begin{equation}
\max \limits_{n\leq N}\|u_n-u_n^k\|_{X} \leq C_4k^2,
\label{tem-err-impli-expli}
\end{equation}
where $C_4>0$ is a constant independent of $k$.
\end{Theorem}
{\it Proof.}
Let $v=u_n$ in \eqref{tem-impli-expli} to get
\begin{equation}
\begin{aligned}
& \langle Au_{n}^{k},u_n-u_{n}^{k} \rangle+\varphi(\cS_{n,L}^{k}u^{k},u_{n}^{k},u_n)-\varphi(\cS_{n,L}^{k}u^{k},u_{n}^{k},u_{n}^{k})\\
& \qquad\quad{} +j^0(\gamma_{j}u_{n}^{k};\gamma_{j}u_n-\gamma_{j}u_{n}^{k})\geq \langle f_n,u_n-u_{n}^{k}\rangle.
\end{aligned}
\label{tem-impli-expli-1}
\end{equation}
Add \eqref{tem-oneorder-3} to \eqref{tem-impli-expli-1} and employ the strong monotonicity of $A$,
\begin{equation}
\begin{aligned}
m_A\|u_n-u_n^k\|_X^2&\le \varphi(\cS_{n,L}^{k}u^{k},u_{n}^{k},u_n)-\varphi(\cS_{n,L}^{k}u^k,u_n^k,u_{n}^{k})\\
&\quad{} +\varphi(\cS_{n}u,u_n,u_n^k)-\varphi(\cS_{n}u,u_n,u_n)\\
&\quad{} +j^0(\gamma_{j}u_n^k;\gamma_{j}u_n-\gamma_{j}u_n^k)+j^0(\gamma_{j}u_n;\gamma_{j}u_n^k-\gamma_{j}u_n)\\
&\le \alpha_\varphi\|u_n-u_n^k\|_X^2+\beta_\varphi\|\cS_nu-\cS_{n,L}^ku^k\|_Y\|u_n-u_n^k\|_X\\
&\quad {}+\alpha_jc_j^2\|u_n-u_n^k\|_X^2.
	\end{aligned}
	\label{tem-iter-eachstep}
\end{equation}
Similar to \eqref{tem-oneorder-6} by using \eqref{eq:S accurancy} instead,
\begin{equation} \label{eq: S two-term}
	\|\cS_nu-\cS_{n,L}^ku^k\|_Y \le Ck^2\|u\|_{W^{2,\infty}(I;X)}+\frac{3}{2}k\beta_{\varphi}\|R\|\|q\|\sum_{j=0}^{n-1}\|u_j-u_j^k\|_X.
\end{equation}
Apply \eqref{eq: S two-term}  to  \eqref{tem-iter-eachstep},
\begin{equation}
\begin{aligned}
\|u_n-u_n^k\|_X&\le  \frac{\beta_\varphi}{m_A-\alpha_\varphi-\alpha_jc_j^2}\|\cS_nu-\cS_{n,L}^ku^k\|_Y \\
& \leq Ck^2\|u\|_{W^{2,\infty}(I;X)}+ \frac{\frac{3}{2}k\beta_{\varphi}\|R\|\|q\|}{m_A-\alpha_{\varphi}-\alpha_{j}c_j^2}\sum_{j=0}^{n-1}\|u_j-u_j^k\|_X.
\end{aligned}
\end{equation}
Then by Corollary \ref{Col Gronwall},
\begin{equation}
\begin{aligned}
\|u_n-u_n^k\|_X&\le \big(k^2C\|u\|_{W^{2,\infty}(I;X)}+\frac{\frac{3}{2}k\beta_{\varphi}\|R\|\|q\|}{m_A-\alpha_{\varphi}-\alpha_{j}c_j^2}\|u_0-u_0^k\|_X\big)\\
&\quad \cdot\bigg(1+k\frac{\frac32\beta_{\varphi}\|R\|\|q\|}{m_A-\alpha_{\varphi}-\alpha_{j}c_j^2}\bigg)^{n-1} \\
& \le C_4k^2,
\end{aligned}
\end{equation}
where
$$ C_4 = C\|u\|_{W^{2,\infty}(I,X)} \cdot \exp\bigg\{ \frac{\frac{3}{2}\beta_{\varphi}\|R\|\|q\|}{m_A-\alpha_{\varphi}-\alpha_{j}c_j^2} t_n\bigg\}.$$
Thus the second-order error estimate \eqref{tem-err-impli-expli} is established.
\hfill$\Box$

\begin{Remark} \label{remark:2rdFix}
For the numerical scheme in \cite{XHHCW19a}, the history dependent operator is implicitly treated in the sense
that its approximation depends on the current unknown solution component.  As a result, a restriction for
the time step-size of the form $k<(m_A-\alpha_{\varphi}-\alpha_{j}c_j^2)/\beta\|R\|\|q\|$ is needed
to ensure the unique solvability and  for the derivation of the error bound there. In contrast,
for our numerical scheme given by Problem \ref{problem tem-impli-expli}, we have the unique solvability and
error bound for an arbitrary time step-size.
\end{Remark}

\smallskip
Next we modify \eqref{tem-3} and give another scheme of second-order.

\begin{Problem} \label{Problem 4.1}
Find a discrete solution $u^{k}:=\{u_{n}^{k}\}_{n=0}^{N}\subset K$ such that
	\begin{equation}
	\begin{aligned}
	& \langle Au_{n}^{k},v-u_{n}^{k} \rangle+\varphi(\cS_{n,L}^{k}u^{k},2u_{n-1}^{k}-u_{n-2}^k,v)-\varphi(\cS_{n,L}^{k}u^{k},2u_{n-1}^{k}-u_{n-2}^k,u_{n}^{k})\\
	& \qquad \quad +j^0(\gamma_{j}u_{n}^{k};\gamma_{j}v-\gamma_{j}u_{n}^{k})+\langle j_c(\gamma_ju_n^k),\gamma_jv-\gamma_ju_n^k\rangle_{X_j} \\
	&\qquad \geq \langle j_c(2\gamma_ju_{n-1}^k-\gamma_ju_{n-2}^k),\gamma_jv-\gamma_ju_n^k\rangle_{X_j}+\langle f_n,v-u_{n}^{k}\rangle \quad \forall\, v \in K, n\ge 2,
	\end{aligned}
	\label{tem-3.1}
	\end{equation}
and for $n=1$,
\begin{equation}
	\begin{aligned}
	& \langle Au_{1}^{k},v-u_{1}^{k} \rangle+\varphi(\cS_{1,L}^{k}u^{k},u_{1}^k,v)-\varphi(\cS_{1,L}^{k}u^{k},u_{1}^k,u_{1}^{k}) \\ &\quad+j^0(\gamma_{j}u_{1}^{k};\gamma_{j}v-\gamma_{j}u_{1}^{k})\geq \langle f_1,v-u_{1}^{k}\rangle \quad \forall\, v \in K.
	\end{aligned}
	\label{tem-onestep-implicit}
\end{equation}
\end{Problem}

The uniqueness and existence results for \eqref{tem-3.1} are similar to that of Problem \ref{Problem 4}.
As for \eqref{tem-onestep-implicit}, it can be referred to Problem \ref{problem tem-impli-expli}.
Then we have the following uniqueness and existence results for Problem \ref{Problem 4.1}.

\begin{Theorem} \label{Theorem uni-exi-2rd-semi-dis}
Assume \eqref{pre-3}--\eqref{pre-10} and \eqref{pre-jc}.  Then Problem \ref{Problem 4.1} has a unique solution $u^k=\{u_n^k\}_{n=0}^N\subset K$.
\end{Theorem}

\smallskip
Next we derive an error bound for semi-discrete solution of Problem~\ref{Problem 4.1}. Meanwhile, a stronger constraint compared with \eqref{eq: small condition 1rd} is needed, i.e.,
\begin{equation}
\alpha_\varphi + \alpha_cc_j^2  < m_A/3. \label{pre-9.1}
\end{equation}

\begin{Theorem} \label{Theorem 7}
Assume \eqref{pre-3}--\eqref{pre-9}, \eqref{pre-jc}, \eqref{pre-9.1} and the regularity
$q\in C^2(\mathbb R_{+} \times \mathbb R_{+};\mathcal{L}(X))$, $u\in {W_{loc}^{2,\infty}(\mathbb R_{+};X)}$.
Then for the semi-discrete solution of Problem \ref{Problem 4.1}, the following error bound holds:
\begin{equation}
  \max \limits_{n\leq N}\|u_n-u_n^k\|_{X} \leq C_5k^2,
  \label{tem-17}
\end{equation}
where $C_5>0$ is a constant independent of $k$.
\end{Theorem}
{\it Proof.}
For $n=1$, we have a second-order accuracy result for \eqref{tem-onestep-implicit} by Theorem \ref{Theorem tem-impli-expli}:
\begin{equation}
\begin{aligned}
\|u_1-u_1^k\|_X \leq C_4k^2.
\end{aligned}
\label{tem-onestep-2nd}
\end{equation}
For $n\ge 2$, taking $v=u_n$ in~\eqref{tem-3.1}, we have
\begin{equation}
\begin{aligned}
& \langle A u_n^k, u_n-u_n^k \rangle + \varphi(\cS_{n,L}^k u^k, 2u_{n-1}^{k}-u_{n-2}^k, u_n)- \varphi(\cS_{n,L}^k u^k, 2u_{n-1}^{k}-u_{n-2}^k, u_n^k)\\
& \qquad\quad {}+j^0(\gamma_{j}u_n^k;\gamma_{j}u_n-\gamma_{j}u_n^k)+\langle j_c(\gamma_ju_n^k),\gamma_ju_n-\gamma_ju_n^k\rangle_{X_j}\\
   &\qquad \geq \langle f_n,u_n-u_n^k \rangle+\langle j_c(2\gamma_ju_{n-1}^{k}-\gamma_ju_{n-2}^k),\gamma_ju_n-\gamma_ju_n^k\rangle_{X_j}.
  \end{aligned}
  \label{tem-20}
\end{equation}
Combine \eqref{tem-oneorder-3} with \eqref{tem-20} and use the strong monotonicity of $A$ to obtain
\begin{equation}
  \begin{aligned}
&m_A\|u_n-u_n^k\|^2_X   \leq \varphi(\cS_{n}u,u_n,u_n^k)-\varphi(\cS_{n}u,u_{n},u_n)\\
 &\qquad\quad+\varphi(\cS_{n,L}^k u^k,2u_{n-1}^{k}-u_{n-2}^k,u_n) -\varphi(\cS_{n,L}^k u^k,2u_{n-1}^{k}-u_{n-2}^k,u_n^k)\\
 &\qquad\quad+j^0(\gamma_{j}u_n;\gamma_{j}u_n^k-\gamma_{j}u_n) +j^0(\gamma_{j}u_n^k;\gamma_{j}u_n-\gamma_{j}u_n^k)\\
 & \qquad\quad+\langle j_c(\gamma_ju_n^k),\gamma_ju_n-\gamma_ju_n^k\rangle_{X_j} -\langle j_c(2\gamma_ju_{n-1}^{k}-\gamma_ju_{n-2}^k),\gamma_ju_n-\gamma_ju_n^k\rangle_{X_j}\\
 &\qquad = \hat E_{\varphi}+E_j+\hat E_{j_c},
  \end{aligned}
  \label{tem-21}
\end{equation}
where $E_j$ is defined in \eqref{tem-Ej-express} and
\begin{align}
\hat E_{\varphi}& =\varphi(\cS_{n,L}^{k}u^{k},2u_{n-1}^{k}-u_{n-2}^{k},u_n)-\varphi(\cS_{n,L}^{k}u^{k},2u_{n-1}^{k}-u_{n-2}^{k},u_{n}^{k})\nonumber \\
& \quad{}+\varphi(\cS_{n}u,u_n,u_{n}^{k})-\varphi(\cS_{n}u,u_n,u_n),
\label{tem-hatEphi-express}\\[1mm]
\hat E_{j_c} &= \langle j_c(\gamma_ju_n),\gamma_ju_n-\gamma_ju_n^k\rangle_{X_j}
-\langle j_c(2\gamma_ju_{n-1}^k-\gamma_ju_{n-2}^k),\gamma_ju_n-\gamma_ju_n^k\rangle_{X_j}.
\label{tem-hatEjc-express}
\end{align}
We bound $\hat E_{j_c}$ and $\hat E_{\varphi}$ as follows:
\begin{equation}
\begin{aligned}
\hat E_{j_c}&=\langle j_c(\gamma_ju_n-2\gamma_ju_{n-1}+\gamma_ju_{n-2}),\gamma_ju_n-\gamma_ju_n^k\rangle_{X_j}\\
&\quad+2\langle j_c(\gamma_ju_{n-1}-\gamma_ju_{n-1}^k),\gamma_ju_n-\gamma_ju_n^k\rangle_{X_j}\\
& \quad-\langle j_c(\gamma_ju_{n-2} -\gamma_ju_{n-2}^k),\gamma_ju_n-\gamma_ju_n^k\rangle_{X_j}\\
&\le \alpha_c c_j^2\big(k^2\|u\|_{W^{2,\infty}(I,X)}+2\|u_{n-1}-u_{n-1}^k\|_X\\
& \quad +\|u_{n-2}-u_{n-2}^k\|_X\big)\|u_n-u_n^k\|_X,
\end{aligned}
\label{tem-hatEjc}
\end{equation}
\begin{equation}
\begin{aligned}
\hat E_{\varphi} &\le \big(k^2\alpha_\varphi\|u\|_{W^{2,\infty}(I,X)}+2\alpha_\varphi \|u_{n-1}-u_{n-1}^k\|_X\\
&\qquad +\alpha_\varphi \|u_{n-2}-u_{n-2}^k\|_X+\beta_\varphi\|\cS_{n}u-\cS_{n,L}^k u^k\|_Y\big) \|u_n-u_n^k\|_X.
\end{aligned}
\label{tem-hatEphi}
\end{equation}
From \eqref{tem-21}--\eqref{tem-hatEphi} and \eqref{eq: S two-term}, we have
\begin{equation}
	\begin{aligned}
m_A\|u_n-u_n^k\|_X&\le Ck^2\|u\|_{W^{2,\infty}(I;X)}+\frac32k\beta_{\varphi}\|R\|\|q\|\sum_{j=0}^{n-1}\|u_j-u_j^k\|_X \\
& \quad +(\alpha_{\varphi}+\alpha_cc_j^2)\big(2\|u_{n-1}-u_{n-1}^k\|_X+\|u_{n-2}-u_{n-2}^k\|_X\big).
	\end{aligned}
	\label{tem-24}
\end{equation}
Apply Lemma \ref{Lemma Gronwall} to \eqref{tem-24} and combine with \eqref{tem-onestep-2nd},
\begin{equation}
\begin{aligned}
\|u_n-u_n^k\|_X&\le \bigg(\frac{\alpha_{\varphi}+\alpha_{c}c_j^2}{m_A}(2\|u_1-u_1^k\|_X+\|u_0-u_0^k\|_X)\\
&\qquad +\frac{\frac32k\beta_{\varphi}\|R\|\|q\|}{m_A-3(\alpha_{\varphi}+\alpha_{c}c_j^2)}\big(\|u_0-u_0^k\|_X+\|u_1-u_1^k\|_X\big) \\
& \qquad+ k^2\frac{C\|u\|_{W^{2,\infty}(I;X)}}{m_A-3(\alpha_\varphi+\alpha_jc_j^2)}\bigg)\cdot\bigg(1+\frac{\frac32k\beta_{\varphi}\|R\|\|q\|}{m_A-3(\alpha_{\varphi}+\alpha_{c}c_j^2)}\bigg)^{n-2} \\
& \le C_5k^2,
\end{aligned}
\end{equation}
where
\begin{align*}
C_5 =& \left(\frac{C\|u\|_{W^{2,\infty}(I;X)}}{m_A-3(\alpha_\varphi+\alpha_cc_j^2)} +  C_4\frac{\frac32k\beta_{\varphi}\|R\|\|q\|}{m_A-3(\alpha_{\varphi}+\alpha_{c}c_j^2)}  \right. \\
& \left. \quad + 2C_4\frac{\alpha_{\varphi}+\alpha_{c}c_j^2}{m_A} \right)  \cdot \exp\bigg\{ \frac{\frac32\beta_{\varphi}\|R\|\|q\|}{m_A-3(\alpha_{\varphi}+\alpha_{c}c_j^2)}t_n \bigg\},
\end{align*}
 which leads to the error bound~\eqref{tem-17}.
\hfill$\Box$

\section{Fully Discrete Approximation}  \label{sec-4}
\setcounter{equation}0

In this section we consider fully discrete approximations of Problem~\ref{Problem 1} with or without constraints. The notation and assumptions follow from previous section, and a regular family of
finite element partitions $\{T^h\} $ with mesh grid size $h$ is introduced for the spatial discretization. Let $X^h \subset X$ be the conforming finite element spaces. We consider internal approximations only, i.e., $K^h=X^h\cap K$ is nonempty, convex and closed.


Certainly, different fully discrete schemes can be constructed with different temporally semi-discrete schemes proposed in the previous section.  We state these fully discrete schemes as follows.

\begin{Problem} \label{Problem ful-1st}
	Find the discrete solution
	$u^{kh}:=\{u_{n}^{kh}\}_{n=0}^{N}\subset K^h$ such that
	\begin{equation}
	\begin{aligned}
	& \langle Au_{n}^{kh},v^h-u_{n}^{kh} \rangle+\varphi(\cS_{n,L}^{k}u^{kh},u_{n-1}^{kh},v^h)-\varphi(\cS_{n,L}^{k}u^{kh},u_{n-1}^{kh},u_{n}^{kh})\\
	&\qquad\quad+j^0(\gamma_{j}u_{n}^{kh};\gamma_{j}v^h-\gamma_{j}u_{n}^{kh})+\langle j_c(\gamma_ju_n^{kh}),\gamma_jv^h-\gamma_ju_n^{kh}\rangle_{X_j} \\
	&\qquad\geq \langle f_n,v^h-u_{n}^{kh}\rangle +\langle j_c(\gamma_ju_{n-1}^{kh}),\gamma_jv^h-\gamma_ju_n^{kh}\rangle_{X_j}\quad \forall\, v^h \in K^h.
	\end{aligned}
	\label{ful-1st}
	\end{equation}
\end{Problem}

\begin{Problem} \label{Problem ful-impli-expli}
	Find the discrete solution
	$u^{kh}:=\{u_{n}^{kh}\}_{n=0}^{N}\subset K^h$ such that
	\begin{equation}
	\begin{aligned}
	& \langle Au_{n}^{kh},v^h-u_{n}^{kh} \rangle+\varphi(\cS_{n,L}^{k}u^{kh},u_{n}^{kh},v^h)-\varphi(\cS_{n,L}^{k}u^{kh},u_{n}^{kh},u_{n}^{kh})\\
	&\qquad+j^0(\gamma_{j}u_{n}^{kh};\gamma_{j}v^h-\gamma_{j}u_{n}^{kh}) \geq \langle f_n,v^h-u_{n}^{kh}\rangle \quad \forall\, v^h \in K^h.
	\end{aligned}
	\label{ful-impli-expli}
	\end{equation}
\end{Problem}

\begin{Problem} \label{Problem 8}
Find the discrete solution $u^{kh}:=\{u_{n}^{kh}\}_{n=0}^N \subset K^h$ such that
\begin{equation}
  \begin{aligned}
   & \langle Au_{n}^{kh},v^h-u_{n}^{kh}\rangle+\varphi(\cS_{n,L}^{k}u^{kh},2u_{n-1}^{kh}-u_{n-2}^{kh},v^h) -\varphi(\cS_{n,L}^{k}u^{kh},2u_{n-1}^{kh}-u_{n-2}^{kh},u_{n}^{kh})\\
   & \qquad\quad +j^0(\gamma_j u_{n}^{kh};\gamma_jv^h-\gamma_ju_{n}^{kh})+\langle j_c(\gamma_ju_n^{kh}),\gamma_jv^h-\gamma_ju_n^{kh}\rangle_{X_j} \\
   &\qquad\ge \langle f_n, v^h-u_{n}^{kh}\rangle + \langle j_c(2\gamma_ju_{n-1}^{kh}-\gamma_ju_{n-2}^{kh}),\gamma_jv^h-\gamma_ju_n^{kh}\rangle_{X_j} \quad\forall\,v^h\in K^h, n\ge 2,
  \end{aligned} \label{ful-3}
\end{equation}
for $n=0, 1$ the following scheme is used
\begin{equation}
\begin{aligned}
& \langle Au_{n}^{kh},v^h-u_{n}^{kh} \rangle+\varphi(\cS_{n,L}^{k}u^{kh},u_{n}^{kh},v^h)-\varphi(\cS_{n,L}^{k}u^{kh},u_{n}^{kh},u_{n}^{kh}) \\ &\quad+j^0(\gamma_{j}u_{n}^{kh};\gamma_{j}v^h-\gamma_{j}u_{n}^{kh})\geq \langle f_n,v^h-u_{n}^{kh}\rangle \quad \forall\, v^h \in K^h.
\end{aligned}
\label{ful-onestep-implicit}
\end{equation}
\end{Problem}

In the following, we will only discuss about the fully discrete Problem \ref{Problem 8}, since
the other two fully discrete schemes can be discussed similarly.
Similar to the temporally semi-discrete case, we can show that under those same
conditions for the temporally semi-discrete, Problem \ref{Problem 8} has a unique solution.
An error bound for Problem~\ref{Problem 8} is given next.

\begin{Theorem} \label{Theorem 9}
Assume \eqref{pre-3}--\eqref{pre-10}, \eqref{pre-jc} and \eqref{pre-9.1}.
Under the regularity assumptions $q\in C^2(\mathbb R_{+} \times \mathbb R_{+};\mathcal{L}(X))$,
$u\in {W_{loc}^{2,\infty}(\mathbb R_{+};X)}$, we have the error bound
\begin{equation}
  \begin{aligned}
    \max\limits_{0\leq n\leq N}\|u_n-u_{n}^{kh} \|_{X} & \leq
   C_6\max\limits_{0\leq n\leq N}\inf \limits_{v^h \in K^h}\{\|u_n-v^h\|_X
   + \|\gamma_ju_{n}-\gamma_jv^h\|_{X_j}^{\frac{1}{2}}\\
  & \quad +|E(v^h,u_n)|^{\frac{1}{2}} \}+ C_6k^2,
  \end{aligned}
  \label{ful-4}
\end{equation}
where $C_6>0$ is a constant independent of $k,h$ and
\begin{equation}
  \begin{aligned}
    E(v^h,u_n) & = \langle Au_n, v^h-u_n \rangle + \varphi(\cS_{n}u,u_n,v^h)- \varphi(\cS_{n}u,u_n,u_n)\\
 & \quad +j^0(\gamma_{j}u_{n};\gamma_{j}v^h-\gamma_{j}u_{n})-\langle f_n, v^h-u_n \rangle, \quad  v^h \in K^h.
  \end{aligned}
  \label{ful-5}
\end{equation}
\end{Theorem}
{\it Proof.}
First we consider the general case of $n\ge 2$. To this end, we take $t=t_n$ and $v= u_{n}^{kh}$
in~\eqref{pre-2} to get
\begin{equation}
  \begin{aligned}
   & \langle Au_n,u_{n}^{kh}-u_n \rangle + \varphi(\cS_{n}u,u_n,u_{n}^{kh})- \varphi(\cS_{n}u,u_n,u_n)\\
   & \quad +j^0(\gamma_{j} u_n;\gamma_{j} u_{n}^{kh}-\gamma_{j}u_n) \geq \langle f_n, u_{n}^{kh}-u_n \rangle.
  \end{aligned}
  \label{ful-6}
\end{equation}
On the other hand,
\begin{equation}
  \begin{aligned}
    \langle Au_n-Au_{n}^{kh}, u_n-u_{n}^{kh} \rangle
    & = \langle Au_n, u_n-u_{n}^{kh} \rangle + \langle Au_{n}^{kh}, u_{n}^{kh}-v^h \rangle \\
    & \quad + \langle Au_{n}^{kh}, v^h-u_n \rangle.
  \end{aligned}
  \label{ful-8}
\end{equation}
Combine \eqref{pre-5}(b) with \eqref{ful-3}, \eqref{ful-6}--\eqref{ful-8},
\begin{equation}
  \begin{aligned}
   m_A\|u_n-u_n^{kh}\|_X^2
   &\leq \varphi(\cS_{n}u,u_n,u_{n}^{kh})-\varphi(\cS_{n}u,u_n,u_n)\\
   & \quad +\varphi(\cS_{n,L}^{k}u^{kh},2u_{n-1}^{kh}-u_{n-2}^{kh},v^h)-\varphi(\cS_{n,L}^{k}u^{kh},2u_{n-1}^{kh}-u_{n-2}^{kh},u_{n}^{kh})\\
   & \quad +\langle Au_{n}^{kh},v^h-u_{n}\rangle-\langle f_n,v^h-u_{n}\rangle+j^0(\gamma_{j}u_n;\gamma_{j}u_{n}^{kh}-\gamma_{j}u_n)\\
   &\quad +\langle j_c(\gamma_ju_n^{kh}),\gamma_jv^h-\gamma_ju_n^{kh}\rangle_{X_j} +j^0(\gamma_{j}u_{n}^{kh};\gamma_{j}v^h-\gamma_{j}u_{n}^{kh})\\
   & \quad -\langle j_c(2\gamma_ju_{n-1}^{kh}-\gamma_ju_{n-2}^{kh}),\gamma_jv^h-\gamma_ju_n^{kh}\rangle_{X_j} \\
   & = E_{\varphi_1}+E_{\varphi_2}+\hat{E}_{j}+E_{A}+E(v^h,u_n),
  \end{aligned}
  \label{ful-9}
\end{equation}
where
\begin{eqnarray}
&& \begin{array}{ll}
    E_{\varphi_1}=\varphi(\cS_{n,L}^{k}u^{kh},2u_{n-1}^{kh}-u_{n-2}^{kh},u_n)-\varphi(\cS_{n,L}^{k}u^{kh},2u_{n-1}^{kh}-u_{n-2}^{kh},u_{n}^{kh}) \\
   \qquad \quad+  \varphi(\cS_{n}u,u_n,u_{n}^{kh})-\varphi(\cS_{n}u,u_n,u_n),
   \end{array}
   \label{ful-10}\\[1mm]
&& \begin{array}{ll}
    E_{\varphi_2}= \varphi(\cS_{n,L}^{k}u^{kh},2u_{n-1}^{kh}-u_{n-2}^{kh},v^h)- \varphi(\cS_{n,L}^{k}u^{kh},2u_{n-1}^{kh}-u_{n-2}^{kh},u_{n}) \\
    \qquad \quad  +\varphi(\cS_{n}u,u_n,u_n) -\varphi(\cS_{n}u,u_n,v^h),
   \end{array}
   \label{ful-11} \\[1mm]
&& \begin{array}{ll}
    \hat{E}_{j} = j^0(\gamma_{j}u_{n};\gamma_{j}u_{n}^{kh}-\gamma_{j}u_{n})+j^0(\gamma_{j}u_{n}^{kh};\gamma_{j}v^h-\gamma_{j}u_{n}^{kh})\\
    \qquad \quad +\langle j_c(\gamma_ju_n),\gamma_ju_n^{kh}-\gamma_ju_n\rangle_{X_j}+\langle j_c(\gamma_ju_n^{kh}),\gamma_ju_n-\gamma_ju_n^{kh}\rangle_{X_j} \\
    \qquad \quad  -j^0(\gamma_{j}u_{n};\gamma_{j}v^h-\gamma_{j}u_{n}),
   \end{array}
   \label{ful-12} \\[1mm]
&& \begin{array}{ll}
	\ E_A = \langle Au_{n}^{kh},v^h-u_n\rangle-\langle Au_n,v^h-u_n \rangle \\ \quad\qquad+\langle j_c(\gamma_ju_n^{kh}),\gamma_jv^h-\gamma_ju_n\rangle_{X_j}-\langle j_c(\gamma_ju_n),\gamma_ju_n^{kh}-\gamma_ju_n\rangle_{X_j} \\
	\qquad \quad -\langle j_c(2\gamma_ju_{n-1}^{kh}-\gamma_ju_{n-2}^{kh}),\gamma_jv^h-\gamma_ju_n^{kh}\rangle_{X_j}.
	\end{array}
   \label{ful-13}
\end{eqnarray}
Let us bound $E_{\varphi_1}$, $E_{\varphi_2}$, $\hat{E}_{j}$ and $E_{A}$ in turn.
\begin{equation}
\begin{aligned}
  E_{\varphi_1}&\leq \alpha_\varphi\|u_n-u_{n}^{kh}\|_X\|u_n-2u_{n-1}^{kh}+u_{n-2}^{kh}\|_X  \\
  & \quad + \beta_\varphi\|\cS_{n}u-\cS_{n,L}^{kh}u^{kh}\|_Y \|u_n-u_{n}^{kh}\|_X,
  \end{aligned}
  \label{ful-14}
\end{equation}
\begin{equation}
  \begin{aligned}
  E_{\varphi_2}&\leq \alpha_\varphi\|u_n-v^h\|_X\|u_n-2u_{n-1}^{kh}+u_{n-2}^{kh}\|_X \\
  & \quad+ \beta_\varphi\|\cS_{n}u-\cS_{n,L}^{kh}u^{kh}\|_Y \|u_n-v^h\|_X.
  \end{aligned}
  \label{ful-15}
\end{equation}
Use the sub-additive property of generalized directional derivative,
\begin{equation}
  \begin{aligned}
 \hat{E}_{j} & \leq j^0(\gamma_{j}u_{n};\gamma_{j}u_n^{kh}-\gamma_{j}v^{h})+j^0(\gamma_{j}u_{n}^{kh};\gamma_{j}v^{h}-\gamma_{j}u_{n}^{kh})-\alpha_jc_j^2\|u_n-u_n^{kh}\|_X^2\\
       & \leq j^0(\gamma_{j}u_n;\gamma_{j}u_n-\gamma_{j}v^h)+j^0(\gamma_{j}u_n;\gamma_{j}u_{n}^{kh}-\gamma_{j}u_n)\\
       & \quad +j^0(\gamma_{j}u_{n}^{kh};\gamma_{j}v^h-\gamma_{j}u_n)+j^0(\gamma_{j}u_{n}^{kh};\gamma_{j}u_n-\gamma_{j}u_{n}^{kh})-\alpha_jc_j^2\|u_n-u_n^{kh}\|_X^2\\
       & \leq (2c_0+c_1\|\gamma_{j}u_n\|_{X_j}+c_1\|\gamma_{j}u_{n}^{kh}\|_{X_j})\|\gamma_{j}u_n-\gamma_{j}v^h\|_{X_j}\\
       & \leq (2c_0+2c_1 c_j\|u_n\|_X)\|\gamma_{j}u_n-\gamma_{j}v^h\|_{X_j}+ c_1 c_j^{2}\|u_n-u_{n}^{kh}\|_X \|u_n-v^h\|_X.
  \end{aligned}
  \label{ful-16}
\end{equation}
Since $A$ is Lipschitz continuous with a Lipschitz constant $L_A>0$,
\begin{equation}
\begin{aligned}
& E_A \leq L_A \|u_n-u_{n}^{kh}\|_X \|u_n-v^h\|_X +\langle j_c(\gamma_ju_n-2\gamma_ju_{n-1}^{kh}+\gamma_ju_{n-2}^{kh}),\gamma_ju_n-\gamma_ju_n^{kh}\rangle_{X_j} \\
&\qquad +\langle j_c(2\gamma_ju_{n-1}^{kh}-\gamma_ju_{n-2}^{kh}),\gamma_ju_n-\gamma_jv^h\rangle_{X_j}+\langle j_c(\gamma_ju_n^{kh}),\gamma_jv^{h}-\gamma_ju_n\rangle_{X_j}.
\end{aligned}
  \label{ful-17}
\end{equation}
We have
\begin{equation}
	\begin{aligned}
	&\langle j_c(2\gamma_ju_{n-1}^{kh}-\gamma_ju_{n-2}^{kh}),\gamma_ju_n-\gamma_jv^h\rangle_{X_j}+\langle j_c(\gamma_ju_n^{kh}),\gamma_jv^{h}-\gamma_ju_n\rangle_{X_j} \\
	&\quad =\langle j_c(-\gamma_ju_n+2\gamma_ju_{n-1}-\gamma_ju_{n-2}),\gamma_ju_n-\gamma_jv^h\rangle_{X_j}\\
	&\quad\quad+\langle j_c(\gamma_ju_n-\gamma_ju_n^{kh}),\gamma_ju_n-\gamma_jv^h\rangle_{X_j} -2\langle j_c(\gamma_ju_{n-1}-\gamma_ju_{n-1}^{kh}),\gamma_ju_n-\gamma_jv^h\rangle_{X_j}\\
	&\quad\quad+\langle j_c(\gamma_ju_{n-2}-\gamma_ju_{n-2}^{kh}),\gamma_ju_n-\gamma_jv^h\rangle_{X_j} \\
	&\quad \le \alpha_cc_j^2\|u_n-v^h\|_X\big(\|u_n-2u_{n-1}+u_{n-2}\|_X+\|u_n-u_n^{kh}\|_X\big) \\
	&\quad \quad +\alpha_cc_j^2\|u_n-v^h\|_X\big(2\|u_{n-1}-u_{n-1}^{kh}\|_X+\|u_{n-2}-u_{n-2}^{kh}\|_X\big).
	\end{aligned}
	\label{ful-18}
\end{equation}
Together with \eqref{ful-9}, \eqref{ful-14}--\eqref{ful-18}, for $\varepsilon< m_A/3 - \alpha_{\varphi} - \alpha_{c}c_j^2$, we obtain
\begin{equation}
	\begin{aligned}
	m_A\|u_n-u_{n}^{kh}\|_X^2
	& \leq \bigg(Ck^2\|u\|_{W^{2,\infty}(I;X)}+2(\alpha_{\varphi}+\alpha_cc_j^2)\|u_{n-1}-u_{n-1}^{kh}\|_X\\
	&+(\alpha_{\varphi}+\alpha_cc_j^2)\|u_{n-2}-u_{n-2}^{kh}\|_X+\beta_{\varphi}\|\cS_{n}u-\cS_{n,L}^{k}u^{kh}\|_Y\\
	&+C\|u_n-v^h\|_X\bigg)\|u_n-u_n^{kh}\|_X +Ck^2\|u_n-v^h\|_X\\
	&+\left( (\alpha_{\varphi}+\alpha_{c}c_j^2)\left(2\|u_{n-1}-u_{n-1}^{kh}\|_X+\|u_{n-2}-u_{n-2}^{kh}\|_X\right) \right. \\
	&  \left. +\beta_{\varphi}\|\cS_{n}u-\cS_{n,L}^{kh}u^{kh}\|_Y^2  \right) \|u_n - v^h\|\\
	&+C\|\gamma_ju_n-\gamma_jv^h\|_{X_j}+|E(v^h,u_n)|.
	\end{aligned}
	\label{ful-19-1}
\end{equation}
Apply Lemma \ref{Lem: square} and Cauchy-Schwarz inequality,
\begin{equation}
  \begin{aligned}
   \|u_n-u_{n}^{kh}\|_X
    & \leq Ck^2\|u\|_{W^{2,\infty}(I;X)}+2(\alpha_{\varphi}+\alpha_cc_j^2)\|u_{n-1}-u_{n-1}^{kh}\|_X\\
    &+(\alpha_{\varphi}+\alpha_cc_j^2)\|u_{n-2}-u_{n-2}^{kh}\|_X+\frac32k\beta_{\varphi}\|R\|\|q\|\sum_{j=0}^{n-1}\|u_j-u_{j}^{kh}\|_X\\
    &+C\|u_n-v^h\|_X + Ck^2 +\frac{\alpha_{\varphi}+\alpha_c^2c_j^2}{m_A\varepsilon}\|u_n-v^h\|_X\\ &+2\varepsilon\|u_{n-1}-u_{n-1}^{kh}\|_X+\varepsilon\|u_{n-2}-u_{n-2}^{kh}\|_X\\
    &+C\|\gamma_ju_n-\gamma_jv^h\|_{X_j}+C|E(v^h,u_n)|.
  \end{aligned}
  \label{ful-19}
\end{equation}

For $n=0$ and $n=1$, a slight modification based on the proof of Theorem \ref{Theorem tem-impli-expli} and the above arguments give
\begin{equation}
\begin{aligned}
&\|u_0-u_{0}^{kh} \|_X
\leq \frac{C}{m_A-\alpha_{\varphi}-\alpha_{j}c_j^2}\{\|u_0-v^h\|_X + \|\gamma_{j}u_{0}-\gamma_{j}v^h\|_{X_j}^{\frac{1}{2}}+|E(v^h,u_0)|^{\frac{1}{2}} \}.
\end{aligned}
\label{ful-0step-err}
\end{equation}
\begin{equation}
\begin{aligned}
&\|u_1-u_{1}^{kh} \|_X
\leq \frac{C}{m_A-\alpha_{\varphi}-\alpha_{j}c_j^2}\{\|u_1-v^h\|_X + \|\gamma_{j}u_{1}-\gamma_{j}v^h\|_{X_j}^{\frac{1}{2}}\\
&\qquad+|E(v^h,u_1)|^{\frac{1}{2}}+k^2\|u\|_{W^{2,\infty}(I;X)}\}+ \frac32k\beta_{\varphi}\|R\|\|q\|\|u_0-u_{0}^{kh}\|_X.
\end{aligned}
\label{ful-1ststep-err}
\end{equation}
Apply Lemma \ref{Lemma Gronwall} to \eqref{ful-19} and combine \eqref{ful-0step-err}--\eqref{ful-1ststep-err} to get
\begin{equation} \label{eq:ful-err}
\begin{aligned}
\|u_n-u_n^{kh}\|_X&\le \left(C\{\|u_n-v^h\|_X + \|\gamma_{j}u_{n}-\gamma_{j}v^h\|_{X_j}^{\frac{1}{2}}+|E(v^h,u_n)|^{\frac{1}{2}} +k^2\|u\|_{W^{2,\infty}(I;X)}\} \right. \\
& \left. \quad+ \frac{\alpha_{\varphi} + \alpha_cc_j^2+\varepsilon}{m_A}(2\|u_1-u_1^{kh}\|_X+\|u_0-u_0^{kh}\|_X) \right.\\
&\left.  \quad +\frac{\frac32k\beta_{\varphi}\|R\|\|q\|}{m_A-3(\alpha_{\varphi}+\alpha_{c}c_j^2+\varepsilon)}\left(\|u_0-u_0^{kh}\|_X+\|u_1-u_1^{kh}\|_X\right)\right)\\
& \left. \qquad\cdot\left(1+\frac{\frac32k\beta_{\varphi}\|R\|\|q\|}{m_A-3(\alpha_{\varphi}+\alpha_{c}c_j^2+\varepsilon)}\right)^{n-2} \right. \\
& \le C_6\max\limits_{0\leq n\leq N}\left(\|u_n-v^h\|_X + \|\gamma_{j}u_{n}-\gamma_{j}v^h\|_{X_j}^{\frac{1}{2}}+|E(v^h,u_n)|^{\frac{1}{2}}+k^2\right),
\end{aligned}
\end{equation}
where
$$ C_6 = C\|u\|_{W^{2,\infty}(I;X)} \cdot \exp\bigg\{ \frac{\frac32\beta_{\varphi}\|R\|\|q\|}{m_A-3(\alpha_{\varphi}+\alpha_{c}c_j^2+\varepsilon)} t_n\bigg\}.$$
Then we have the error bound~\eqref{ful-4}.
\hfill$\Box$

\smallskip

Now we consider the error estimation for numerical solution of the discrete problem without constraint.
We introduce the following assumption on $\varphi$ as in \cite{XHHCW19a}, which allows us to simplify
the error bound~\eqref{ful-4}:
\begin{equation}
  \begin{aligned}
  \left \{
     \begin{array}{ll}
       \varphi:Y\times K\times K \rightarrow \mathbb R  \ \textrm{is a function such that}\\
       \quad  \textrm{there exists a constant} \ c_\varphi >0 \ \textrm{satisfies}\\
       \quad  \varphi(y,u,v_1)+\varphi(y,u,v_2)-2\varphi(y,u,{\frac{v_1+v_2}{2}})\leq c_\varphi\|v_1-v_2\|_X^2 \\
       \qquad \forall\, y \in Y, \ \forall\, u, v_1, v_2 \in K.
     \end{array}
    \right.
  \end{aligned}
  \label{ful-1}
\end{equation}

\begin{Theorem} \label{Theorem 10}
Keep the assumptions stated in Theorem \ref{Theorem 9}. In addition, let $K=X$ and the function
$\varphi$ satisfy the assumption \eqref{ful-1}.  Then the following error bound holds:
\begin{equation}
  \begin{aligned}
    \max \limits_{0\leq n\leq N}\|u_n-u_{n}^{kh}\|_X & \le C\left(\max\limits_{0\leq n\leq N}
		\inf \limits_{v^h\in K^h}\{\|u_n-v^h\|_X + \|\gamma_ju_{n}-\gamma_jv^h\|_{X_j}^{\frac{1}{2}}\}+ k^2\right).\\
  \end{aligned}
  \label{ful-22}
\end{equation}
\end{Theorem}
{\it Proof.}
We start with
\begin{equation}
  \begin{aligned}
  \langle Au_n-Au_{n}^{kh},u_n-u_{n}^{kh}\rangle
  & =\langle Au_n-Au_{n}^{kh},u_n-v^{h}\rangle+\langle Au_n-Au_{n}^{kh},v^{h}-u_{n}^{kh}\rangle\\
  & = \langle Au_n-Au_{n}^{kh},u_n-v^{h}\rangle +\langle Au_n,v^{h}-u_{n}\rangle\\
  & \quad +\langle Au_n,u_{n}-u_{n}^{kh}\rangle+\langle Au_{n}^{kh},u_{n}^{kh}-v^{h}\rangle.\\
  \end{aligned}
  \label{ful-23}
\end{equation}
Further, we replace $v$ with $2u_n-v$ in \eqref{tem-oneorder-1} to get
\begin{equation}
  \begin{aligned}
   & \langle Au_n, u_n-v \rangle + \varphi(\cS_{n}u,u_n,2u_n-v) - \varphi(\cS_{n}u,u_n,u_n)\\
   & \quad +j^0(\gamma_ju_n;\gamma_ju_n-\gamma_{j}v) \geq \langle f_n, u_n-v\rangle \quad \forall\, v \in X.
  \end{aligned}
  \label{ful-25}
\end{equation}
Similarly, take $v=v^{h}$ in~\eqref{ful-25} to get
\begin{equation}
  \begin{aligned}
   & \langle Au_n, u_n-v^h \rangle + \varphi(\cS_{n}u,u_n,2u_n-v^h) - \varphi(\cS_{n}u,u_n,u_n)\\
   & \quad +j^0(\gamma_ju_n; \gamma_ju_n-\gamma_jv^h) \geq \langle f_n, u_n-v^h \rangle.
  \end{aligned}
  \label{ful-26}
\end{equation}
Combine \eqref{pre-5}, \eqref{ful-3}, \eqref{ful-6}, \eqref{ful-23} and \eqref{ful-26},
\begin{equation}
  \begin{aligned}
    m_A\|u_n-u_n^{kh}\|_X^2
    & \leq \langle Au_n-Au_{n}^{kh}, u_n-v^h \rangle + \varphi(\cS_{n}u,u_{n},2u_{n}-v^h)\\
    & \quad + \varphi(\cS_{n}u,u_{n},u_{n}^{kh})-2\varphi(\cS_{n}u,u_{n},u_{n})\\
    & \quad + \varphi(\cS_{n,L}^{k}u^{kh},2u_{n-1}^{kh}-u_{n-2}^{kh},v^h)- \varphi(\cS_{n,L}^{k}u^{kh},2u_{n-1}^{kh}-u_{n-2}^{kh},u_{n}^{kh})\\
    & \quad +j^0(\gamma_ju_{n};\gamma_ju_{n}-\gamma_jv^h)+j^0(\gamma_ju_{n};\gamma_ju_{n}^{kh}-\gamma_ju_{n})\\
    & \quad +j^0(\gamma_ju_{n}^{kh};\gamma_jv^h-\gamma_ju_{n}^{kh})+\langle j_c(\gamma_ju_n^{kh}),\gamma_jv^h-\gamma_ju_n^{kh}\rangle_{X_j}\\
    &\quad -\langle j_c(2\gamma_ju_{n-1}^{kh}-\gamma_ju_{n-2}^{kh}),\gamma_jv^h-\gamma_ju_n^{kh}\rangle_{X_j}\\
    & =E_{\varphi_1}+E_{\varphi_2}+E_{\varphi_3}+\widetilde{E}_{j}+E_{A},
  \end{aligned}
  \label{ful-27}
\end{equation}
where $E_{\varphi_1}$, $E_{\varphi_2}$, $E_{A}$ are the same as in \eqref{ful-10}, \eqref{ful-11}, \eqref{ful-13}
respectively with their bounds \eqref{ful-14}, \eqref{ful-15}, \eqref{ful-17}. In addition,
\begin{equation}
  \begin{aligned}
  E_{\varphi_3} = \varphi(\cS_{n}u,u_n,2u_{n}-v^h) + \varphi(\cS_{n}u,u_{n},v^h) - 2\varphi(\cS_{n}u,u_n,u_n),
  \end{aligned}
  \label{ful-28}
\end{equation}
\begin{equation}
  \begin{aligned}
  \widetilde{E}_{j}
 & = j^0(\gamma_{j}u_{n}; \gamma_{j}u_{n}-\gamma_{j}v^h) +j^0(\gamma_{j}u_{n};\gamma_{j}u_{n}^{kh}-\gamma_{j}u_{n})\\
    & \quad +j^0(\gamma_{j}u_{n}^{kh};\gamma_{j}v^h-\gamma_{j}u_{n}^{kh})-\langle j_c(\gamma_ju_n-\gamma_ju_n^{kh}),\gamma_ju_n-\gamma_ju_n^{kh}\rangle_{X_j}.
  \end{aligned}
  \label{ful-29}
\end{equation}
The assumption~\eqref{ful-1} shows that
\begin{equation}
  \begin{aligned}
   E_{\varphi_3} \leq C \|u_n-v^h\|_X^2.
  \end{aligned}
  \label{ful-30}
\end{equation}
Using the sub-additive property again, we obtain
\begin{equation}
  \begin{aligned}
   \widetilde{E}_{j} \leq  C\|\gamma_{j}u_n-\gamma_{j}v^h\|_{X_j}+C\|u_n-u_{n}^{kh}\|_X\|u_n-v^h\|_X.
  \end{aligned}
  \label{ful-31}
\end{equation}
Together \eqref{ful-27} with \eqref{ful-28}--\eqref{ful-31} and analogy to  \eqref{ful-19},
\begin{equation}
  \begin{aligned}
   m_A \|u_n-u_{n}^{kh}\|_X & \leq C\{\|u_n-v^h\|_X + \|\gamma_ju_{n}-\gamma_jv^h\|_{X_j}^{\frac{1}{2}}\}+ C\|\cS_{n}u-\cS_{n,L}^{k}u^{kh}\|_Y \\
   &\quad +(\alpha_cc_j^2+\alpha_{\varphi}+\varepsilon)\big(\|u_{n-2}-u_{n-2}^{kh}\|_X+2\|u_{n-1}-u_{n-1}^{kh}\|_X\big).
  \end{aligned}
  \label{ful-32}
\end{equation}
Similar to the constrained situation, the error bounds for $n=0,1$ are
\begin{equation}
\|u_0-u_{0}^{kh} \|_X
\leq \frac{C}{m_A-\alpha_{\varphi}-\alpha_{c}c_j^2}\{\|u_0-v^h\|_X + \|\gamma_{j}u_{0}-\gamma_{j}v^h\|_{X_j}^{\frac{1}{2}} \},
\label{ful-0step-err-uncontrain}
\end{equation}
\begin{equation}
\begin{aligned}
\|u_1-u_{1}^{kh} \|_X
&\leq \frac{C}{m_A-\alpha_{\varphi}-\alpha_{c}c_j^2}\{\|u_1-v^h\|_X + \|\gamma_{j}u_{1}-\gamma_{j}v^h\|_{X_j}^{\frac{1}{2}}\\
&\quad+k^2\|u\|_{W^{2,\infty}(I;X)}\}+ Ck\|u_0-u_{0}^{kh}\|_X.
\end{aligned}
\label{ful-1ststep-err-uncontrain}
\end{equation}
Combining  \eqref{ful-32}--\eqref{ful-1ststep-err-uncontrain}, we find
the following error bound by an application of Lemma \ref{Lemma Gronwall},
\begin{equation} \label{eq:ful-err-all}
\begin{aligned}
\|u_n-u_n^{kh}\|_X&\le \left(C\{\|u_n-v^h\|_X + \|\gamma_{j}u_{n}-\gamma_{j}v^h\|_{X_j}^{\frac{1}{2}}+ +k^2\|u\|_{W^{2,\infty}(I;X)}\} \right. \\
& \left. \quad+ \frac{\alpha_{\varphi} + \alpha_cc_j^2+\varepsilon}{m_A}(2\|u_1-u_1^{kh}\|_X+\|u_0-u_0^{kh}\|_X) \right.\\
&\left.  \quad +\frac{Ck}{m_A-3(\alpha_{\varphi}+\alpha_{c}c_j^2+\varepsilon)}\left(\|u_0-u_0^{kh}\|_X+\|u_1-u_1^{kh}\|_X\right)\right)\\
& \left. \qquad\cdot\left(1+\frac{Ck}{m_A-3(\alpha_{\varphi}+\alpha_{c}c_j^2+\varepsilon)}\right)^{n-2} \right. \\
& \le C\max\limits_{0\leq n\leq N}\left(\|u_n-v^h\|_X + \|\gamma_{j}u_{n}-\gamma_{j}v^h\|_{X_j}^{\frac{1}{2}}+k^2\right).
\end{aligned}
\end{equation}
Thus, the proof is completed.
\hfill$\Box$
%

\section{Numerical computation using fixed-point iteration} \label{sec-6}
\setcounter{equation}0

Notice that in Problems \ref{problem tem-impli-expli}, \ref{Problem ful-impli-expli} and in the initial steps of Problems \ref{Problem 4.1}, \ref{Problem 8}, the implicit discretization with respect to the unknown solution
component is used.  Let us discuss how to implement these numerical schemes in practice. We use
a fixed-point iteration approach.  We first consider the fixed-point iterations for the temporally
semi-discrete schemes.

\begin{Problem} \label{Problem 10}
Let $TOL$ be a given error tolerance.  For $1\le n\le N$, find a sequence $\{\tilde{u}_{n,i}^{k}\} \subset K$
from the iterations
\begin{equation}
	\begin{aligned}
	& \langle A\tilde{u}_{n,i}^k,v-\tilde{u}_{n,i}^k \rangle+\varphi(\cS_{n,L}^ku^k,\tilde{u}_{n,i-1}^k,v)-\varphi(\cS_{n,L}^ku^k,\tilde{u}_{n,i-1}^k,\tilde{u}_{n,i}^k)\\
	& \qquad\quad +j^0(\gamma_{j}\tilde{u}_{n,i}^k;\gamma_jv-\gamma_{j}\tilde{u}_{n,i}^k)+\langle j_c(\gamma_j\tilde{u}_{n,i}^k),\gamma_jv-\gamma_j\tilde{u}_{n,i}^k\rangle_{X_j}\\
	&\quad \geq \langle f_n,v- \tilde{u}_{n,i}^k \rangle+\langle j_c(\gamma_j\tilde{u}_{n,i-1}^k),\gamma_jv-\gamma_j\tilde{u}_{n,i}^k\rangle_{X_j},\quad \forall v\in K
	\end{aligned}
	\label{tem-onestep-explicit}
	\end{equation}
until the relative error $\frac{\|\tilde{u}_{n,i}^k-\tilde{u}_{n,i-1}^k\|_X}{\|\tilde{u}_{n,i}^k\|_X}<TOL$;
choose $u^k_{n}$ to be the last iteration $\tilde{u}^k_{n,i}$.
\end{Problem}

\smallskip
In Problem \ref{Problem 10}, the index $i$ refers to the $i$-th iterate at time level $t_n$. For the initialization of iteration, we may use the iterative solution from the previous step, i.e.,
$\tilde{u}_{n,0}^k = u_{n-1}^k$ for $n\ge1$. Now we consider the convergence of the sequence $\{\tilde{u}_{n,i}^{k}\}$ generated by \eqref{tem-onestep-explicit} to the solution of \eqref{tem-impli-expli}.

\begin{Theorem} \label{Therorem-tem-fix-point}
Assume \eqref{pre-3}--\eqref{pre-10}.  Then the iteration \eqref{tem-onestep-explicit} converges linearly with
a convergence rate $\rho=(\alpha_\varphi+\alpha_cc_j^2)/m_A$ that is independent of the time step-size $k$.
\end{Theorem}
{\it Proof.}
Take $v=\tilde{u}_{n,i}^k$ in \eqref{tem-impli-expli},
\begin{equation}
\begin{aligned}
& \langle Au_n^k,\tilde{u}_{n,i}^k-u_n^k \rangle+\varphi(\cS_{n,L}^ku^k,u_n^k,\tilde{u}_{n,i}^k)-\varphi(\cS_{n,L}^ku^k,u_n^k,u_n^k)\\
& \quad{} +j^0(\gamma_{j}u_n^k;\gamma_{j}\tilde{u}_{n,i}^k-\gamma_{j}u_n^k)\ge\langle f_n,\tilde{u}_{n,i}^k - u_n^k \rangle.
\end{aligned}
\label{tem-fixed-iter}
\end{equation}
Take $v=u_n^k$ in \eqref{tem-onestep-explicit},
\begin{equation}
\begin{aligned}
& \langle A\tilde{u}_{n,i}^k,u_n^k-\tilde{u}_{n,i}^k \rangle+\varphi(\cS_{n,L}^ku^k,\tilde{u}_{n,i-1}^k,u_n^k)-\varphi(\cS_{n,L}^ku^k,\tilde{u}_{n,i-1}^k,\tilde{u}_{n,i}^k)\\
	& \qquad\quad +j^0(\gamma_{j}\tilde{u}_{n,i}^k;\gamma_{j}u_n^k-\gamma_{j}\tilde{u}_{n,i}^k)+\langle j_c(\gamma_j\tilde{u}_{n,i}^k),\gamma_ju_n^k-\gamma_j\tilde{u}_{n,i}^k\rangle_{X_j}\\
	&\quad \geq \langle f_n,u_n^k- \tilde{u}_{n,i}^k \rangle+\langle j_c(\gamma_j\tilde{u}_{n,i-1}^k),\gamma_ju_n^k-\gamma_j\tilde{u}_{n,i}^k\rangle_{X_j}.
\end{aligned}
\label{tem-fixed-iter-dis}
\end{equation}
Combine \eqref{tem-fixed-iter} with \eqref{tem-fixed-iter-dis},
\begin{equation}
\begin{aligned}
	& \langle A\tilde{u}_{n,i}^k,u_n^k-\tilde{u}_{n,i}^k \rangle+\langle Au_n^k,\tilde{u}_{n,i}^k-u_n^k \rangle \le \alpha_{\varphi}\|u_n^k-\tilde{u}_{n,i}^k\|_X\|u_n^k-\tilde{u}_{n,i-1}^k\|_X \\
	&\qquad\quad +\alpha_cc_j^2\|u_n-\tilde{u}_{n,i}^k\|_X\|u_n^k-\tilde{u}_{n,i-1}^k\|_X.
	\end{aligned}
	\label{tem-one-iter-app}
	\end{equation}
By the strong monotonicity of $A$ and \eqref{tem-one-iter-app}, we have the following relation:
\begin{equation}
m_A\|u_n^k-\tilde{u}_{n,i}^k\|_X \le \big(\alpha_{\varphi}+\alpha_cc_j^2\big)\|u_n^k-\tilde{u}_{n,i-1}^k\|_X.
\label{tem-one-step-app}
\end{equation}
Therefore, the stated result is proved.
\hfill$\Box$	

	\smallskip
In analogy to the temporally semi-discrete scheme, the iteration algorithm for the fully discrete scheme can be stated as follows.

\begin{Problem} \label{Problem 15}
Let $TOL$ be a given error tolerance.  For $1\le n\le N$, find a sequence $\{\tilde{u}_{n,i}^{kh}\} \subset K^h$ such that
\begin{equation}
\begin{aligned}
& \langle A\tilde{u}_{n,i}^{kh},v^h-\tilde{u}_{n,i}^{kh} \rangle+\varphi(\cS_{n,L}^{k}u^{kh},\tilde{u}_{n,i-1}^{kh},v^h)-\varphi(\cS_{n,L}^{k}u^{kh},\tilde{u}_{n,i-1}^{kh},\tilde{u}_{n,i}^{kh})\\
& \qquad\quad +j^0(\gamma_{j}\tilde{u}_{n,i}^{kh};\gamma_jv^h-\gamma_{j}\tilde{u}_{n,i}^{kh})+\langle j_c(\gamma_j\tilde{u}_{n,i}^{kh}),\gamma_jv^h-\gamma_j\tilde{u}_{n,i}^{kh}\rangle_{X_j}\\
&\quad \geq \langle f_n,v^h- \tilde{u}_{n,i}^{kh} \rangle+\langle j_c(\gamma_j\tilde{u}_{n,i-1}^{kh}),\gamma_jv^h-\gamma_j\tilde{u}_{n,i}^{kh}\rangle_{X_j}, \quad \forall v^h\in K^h,
\end{aligned}
\label{ful-onestep-explicit}
\end{equation}
until the relative error $\frac{\|\tilde{u}_{n,i}^{kh}-\tilde{u}_{n,i-1}^{kh}\|_X}{\|\tilde{u}_{n,i}^{kh}\|_X}<TOL$; choose $u^{kh}_{n}$ to be the last iteration $\tilde{u}^{kh}_{n,i}$.
\end{Problem}
	
The sequence $\{\tilde{u}_{n,i}^{kh}\}$ can be similarly proved to converge to the solution of \eqref{ful-impli-expli}.

\begin{Theorem}\label{Therorem-ful-fix-point}
Keep the assumptions in Theorem \ref{Therorem-tem-fix-point}.  Then the iteration \eqref{ful-onestep-explicit} converges linearly with a convergence rate $\rho=(\alpha_\varphi+\alpha_cc_j^2)/m_A$ that is independent of the time step $k$ and the mesh parameter $h$.
\end{Theorem}

So far we have proposed three types of schemes and the corresponding numerical treatments to solve
Problem \ref{Problem 1}. Note that the difference of the schemes lies in the way the temporal discretization
is done.  We list the schemes
and summarize their main properties in Table \ref{Tab-sum}, where CO stands for convergence order.

\begin{table}[htbp]
	\centering
	\caption{Comparison of the three temporally semi-discrete schemes}
	\begin{tabular}{|c|c|c|c|}
		\hline
		semi-discrete problem & numerical method & CO & constraint  \\
		\hline
		Problem \ref{Problem 4} & \tabincell{l}{$\bm\cdot$ convex optimization} & first-order & $m_A>\alpha_\varphi+\alpha_cc_j^2$ \\
		\hline
		Problem \ref{problem tem-impli-expli} &\tabincell{l}{$\bm\cdot$ convex optimization \\ $\bm\cdot$ fixed-point iteration\\
		\quad (each step)} & second-order & $m_A>\alpha_\varphi+\alpha_jc_j^2$ \\
		\hline
		Problem \ref{Problem 4.1} & \tabincell{l}{$\bm\cdot$ convex optimization \\ $\bm\cdot$ extrapolation\\
		$\bm\cdot$ fixed-point iteration\\
		\quad (initial step)} & second-order & $m_A/3>\alpha_\varphi+\alpha_cc_j^2$ \\
		\hline
	\end{tabular}
	\label{Tab-sum}
\end{table}
We use the result of previous step to approximate the current step in Problem \ref{Problem 4} which is easy to implement while with low accuracy. For Problem \ref{Problem 4.1}, the approximation for current step is performed with an extrapolation, thus an initial step is introduced and we employ a fixed-point iteration to solve it numerically. As a result, we obtain a second-order accuracy with stronger small condition constraint. Inspired by this fixed-point iterative procedure, we propose a new scheme in Problem \ref{problem tem-impli-expli},  in which a fixed-point iteration is used to approximate this scheme for each step.

\section{Application to a contact problem} \label{sec-5}
\setcounter{equation}0

In this section we apply the abstract numerical analysis results in the previous sections to a particular
history-dependent variational-hemivariational inequality.  A viscoelastic frictionless contact model studied
in~\cite{SM16} will be considered.  For details on the model, we refer the reader to \cite{SM16,XHHCW19a}.

\begin{Problem} \label{Problem 11}
Find a displacement $\textbf{u}:\Omega\times \mathbb R_{+}  \rightarrow\mathbb R^d $ and a stress field
$\boldsymbol{\sigma}:\Omega\times \mathbb R_{+}\rightarrow\mathbb S^d$ such that for all $t\in\mathbb R_{+}$,
\begin{eqnarray}
&& \begin{array}{ll}
   \boldsymbol{\sigma}(t)=\mathscr{A}\boldsymbol{\varepsilon}(\boldsymbol{u}(t))
+\mu(\boldsymbol{\varepsilon}(\boldsymbol{u}(t))-P_{M(\kappa(\zeta(t)))}\boldsymbol{\varepsilon}(\boldsymbol{u}(t)) )\\
   \qquad\quad +\int_{0}^{t}\mathscr{B}(t-s)\boldsymbol{\varepsilon}(\boldsymbol{u}(s))ds \quad {\rm in} \ \Omega,
   \end{array}
   \label{err-01}\\[1mm]
&& \begin{array}{ll}
   \mathrm{Div} \ \boldsymbol{\sigma}(t)+\boldsymbol{f_0}(t)=\mathbf{0} \quad {\rm \ in} \ \Omega,
   \end{array}
   \label{err-02} \\[1mm]
&& \begin{array}{ll}
   \textbf{u}(t)=\mathbf{0} \quad {\rm \ on \ } \Gamma_1,
   \end{array}
   \label{err-03} \\[1mm]
&& \boldsymbol{\sigma}(t)\boldsymbol{\nu}=\boldsymbol{f_2}(t) \quad  {\rm \ on \ } \Gamma_2,
   \label{err-04} \\[1mm]
&&  \left \{
    \begin{array}{ll}
    u_\nu(t)\leq g, \ \sigma_\nu(t)+\xi_\nu(t)\leq 0, \\
    (\sigma_\nu(t)+\xi_\nu(t))(u_\nu(t)-g)=0, \quad  {\rm on} \ \Gamma_3, \\
    \xi_\nu(t)\in\partial j_\nu(u_\nu(t)) \\
    \end{array}
    \right.
    \label{err-05} \\[1mm]
&& \begin{array}{ll}
    \boldsymbol{\sigma}_{\boldsymbol{\tau}}(t)=\mathbf{0} \quad {\rm on} \ \Gamma_3.
    \end{array}
    \label{err-06}
\end{eqnarray}
\end{Problem}

\smallskip
As is standard in the literature in the area of the paper, we denote by $\mathbb{S}^d$ the space of second order symmetric tensors on $\mathbb{R}^d$, $\boldsymbol{u}=(u_i)$,
$\boldsymbol{\nu}=(\nu_i)$, $\boldsymbol{\sigma}=(\sigma_{ij})$,
$\boldsymbol{\varepsilon}(\boldsymbol{u})=(\nabla \boldsymbol{u} +(\nabla \boldsymbol{u})^T)/2 $
the displacement field, outward unit normal on the boundary, stress tensor and linearized strain tensor,
respectively. In addition, $v_\nu:=\boldsymbol{v}\cdot \boldsymbol{\nu}$ and
$\boldsymbol{v}_{\tau}:=\boldsymbol{v}- v_\nu \boldsymbol{\nu}$ stand for the normal and
tangential components of a vector field $\boldsymbol{v}$, $\sigma_\nu:=(\boldsymbol{\sigma \nu})
\cdot \boldsymbol{\nu}$ and $\boldsymbol{\sigma}_{\tau}:=\boldsymbol{\sigma \nu}- \sigma_\nu \boldsymbol{\nu}$
represent the normal and tangential components of the stress field $\boldsymbol{\sigma}$, respectively.
In equation~(\ref{err-01}) $ P_{M(\kappa(\cdot))}$ denotes the projection on the Von Mises convex,
$\mathscr{A}$  and $\mathscr{B}$ are the elastic and relaxation tensors, and $\mu$ is a constant.
In this model, time-dependent surface tractions of density $\boldsymbol{f}_2$ and volume forces of density $\boldsymbol{f}_0$ are considered. On $\Gamma_3$, the penetration is restricted by a non-negative function
$g$ and the potential function is denoted as $j_{\nu}$.  The function spaces $V$ and $\mathcal{H}$ are
\begin{align*}
& V=\{\boldsymbol{v}=(v_i)\in H^1(\Omega;\mathbb R^{d})\mid \boldsymbol{v}=0 \ {\rm a.e.\ on} \ \Gamma_1\},\\
& \mathcal{H}=\{\boldsymbol{\tau}=(\tau_{ij})\in L^2(\Omega;\mathbb S^{d})\mid
\tau_{ij}= \tau_{ji}, \ 1\leq i, j\leq d \}.
\end{align*}
The inner products in the Hilbert spaces $\mathcal{H}$ and $V$ are
\[(\boldsymbol{\sigma}, \boldsymbol{\tau})_{\mathcal{H}}=\int_{\Omega}\sigma_{ij}(x)\tau_{ij}(x)dx,\quad (\boldsymbol{u}, \boldsymbol{v})_{V}=(\boldsymbol{\varepsilon}(\boldsymbol{u}),
\boldsymbol{\varepsilon}(\boldsymbol{u}))_{\mathcal{H}} \]
and the associated norm are denoted by $\|\cdot\|_{\mathcal{H}}$
and $\|\cdot\|_{V}$. The space of fourth order tensor fields $\mathrm{Q_\infty}$ is given by
\[\mathrm{Q_\infty}=\{ \mathcal{E}=(\mathcal{E}_{ijkl})\mid
\mathcal{E}_{ijkl}=\mathcal{E}_{jikl}=\mathcal{E}_{klij} \in L^\infty(\Omega), \ 1\leq i,j,k,l \leq d \}.\]

We now list the assumptions on the problem data, following \cite{SM16,XHHCW19a}. The elasticity tensor
$\mathcal{A}: \Omega\times\mathbb{S}^{d} \rightarrow \mathbb{S}^{d}$ is symmetric and positive. The relaxation tensor $\mathcal{B}\in C(\mathbb{R_{+}}; Q_{\infty})$ and the bound $\kappa:\mathbb{R}\rightarrow\mathbb{R_{+}}$
is Lipschitz continuous.  The potential function $j_{\nu}: \Gamma_3\times\mathbb{R} \rightarrow \mathbb{R}$ is measurable with respect to the first argument on $\Gamma_3$ for all $r\in\mathbb{R}$ and is locally Lipschitz
with respect to the second argument on $\mathbb{R}$ for a.e.\ $\bm{x}\in \Gamma_3$;
$j_{\nu}(\cdot, \overline{e}(\cdot))$ belongs to $L^1(\Gamma_3)$ for some $\overline{e}\in L^2(\Gamma_3)$. Besides, $|\partial j_{\nu}(\bm{x}, r)|\le \overline{c}_0+\overline{c}_1|r|$ for a.e.\ $\bm{x}\in \Gamma_3$,
for all $r\in \mathbb{R}$ with $\overline{c}_0, \overline{c}_1>0$. In addition, there exists
$\overline{\alpha}_{\nu} \ge 0$ such that for a.e.\ $\bm{x} \in \Gamma_3$,
$$ j_{\nu}^0(\bm{x}, r_1; r_2-r_1) +j_{\nu}^0(\bm{x}, r_2; r_1-r_2) \le \overline{\alpha}_{\nu}|r_1-r_2|^2 \quad \forall\, r_1,r_2 \in \mathbb{R}.$$
For the body force and surface traction, we assume $\bm{f}_0 \in C(\mathbb{R_{+}}; L^2(\Omega; \mathbb{R}^d))$ and  $\bm{f}_2\in C(\mathbb{R_{+}}; L^2(\Gamma_2; \mathbb{R}^d))$.
Let
$ U=\{ \boldsymbol{v} \in V\mid v_\nu\leq g \ {\rm a.e. \ on } \ \Gamma_3 \}$
be the set of admissible displacements. Define the function $\boldsymbol{f}: \mathbb{R_+} \rightarrow V^*$ by
\[\langle \boldsymbol{f}(t), \boldsymbol{v} \rangle_{V^*\times V}=(\boldsymbol{f_0}(t), \boldsymbol{v} )_{L^2(\Omega;\mathbb R^d)}
 +(\boldsymbol{f_2}(t), \boldsymbol{v})_{L^2(\Gamma_2;\mathbb R^d)} \quad \forall\, \boldsymbol{v}\in V,\
 \forall\, t \in \mathbb{R_{+}}. \]
Then the weak formulation of Problem \ref{Problem 11} can be described as following.

\begin{Problem} \label{Problem 12}
Find a displacement $\textbf{u}:\mathbb R_{+} \rightarrow U $ such that the following inequality holds:
\begin{equation}
  \begin{aligned}
   & (\mathscr{A}\boldsymbol{\varepsilon}(\boldsymbol{u}(t)), \boldsymbol{\varepsilon}(\boldsymbol{v})-\boldsymbol{\varepsilon}(\boldsymbol{u}(t)))_{\mathcal{H}}
   +\mu(\boldsymbol{\varepsilon}(\boldsymbol{u}(t)), \boldsymbol{\varepsilon}(\boldsymbol{v})-\boldsymbol{\varepsilon}(\boldsymbol{u}(t)))_{\mathcal{H}}\\
   & \quad -\mu(P_{M(\kappa(\zeta(t)))}\boldsymbol{\varepsilon}(\boldsymbol{u}(t)), \boldsymbol{\varepsilon}(\boldsymbol{v})-\boldsymbol{\varepsilon}(\boldsymbol{u}(t)) )_{\mathcal{H}} \\
   & \quad +\left(\int_{0}^{t}\mathscr{B}(t-s)\boldsymbol{\varepsilon}(\boldsymbol{u}(s))ds, \boldsymbol{\varepsilon}(\boldsymbol{v})-\boldsymbol{\varepsilon}(\boldsymbol{u}(t)) \right)_{\mathcal{H}} \\
   & \quad +\int_{\Gamma_3} j_v^0(u_{\nu}(t); v_{\nu}-u_{\nu}(t))d\Gamma
     \geq \langle \textbf{f}(t), \textbf{v}-\textbf{u}(t) \rangle_{V^*\times V}
    \quad \forall\, \textbf{v} \in U, \ t \in \mathbb R_{+}. \\
  \end{aligned}
  \label{err-07}
\end{equation}
\end{Problem}

To apply the abstract results from the previous sections to the study of this contact problem, some
definitions are needed. We let $\gamma_{j}:V \rightarrow L^2(\Gamma_3)$ be the trace operator defined
by $\gamma_j\boldsymbol{v}=\boldsymbol{v}_\nu$ for $\boldsymbol{v} \in V$. In addition, we define the following operators (\cite{SM16,XHHCW19a}):
\begin{equation}
  \begin{aligned}
  \langle A \boldsymbol{u},\boldsymbol{v} \rangle_{V^*\times V}
    =(\mathscr{A}\boldsymbol{\varepsilon}(\boldsymbol{u}), \boldsymbol{\varepsilon}(\boldsymbol{v} ))_{\mathcal{H}}
    +\mu(\boldsymbol{\varepsilon}(\boldsymbol{u}), \boldsymbol{\varepsilon}(\boldsymbol{v}))_{\mathcal{H}}\\
  \end{aligned}
  \label{err-8}
  \quad \forall\, \boldsymbol{u},\boldsymbol{v} \in V,
\end{equation}
\begin{equation}
  \begin{aligned}
  \|y\|_Y=|r|+\|\boldsymbol{\theta}\|_{\mathcal{H}}
  \quad \forall\, y=(r, \boldsymbol{\theta})\in Y:=\mathbb{R}\times\mathcal{H},
  \end{aligned}
  \label{err-9}
\end{equation}
\begin{equation}
  \begin{aligned}
  \varphi(y,\boldsymbol{u},\boldsymbol{v})
  &= -\mu(P_{M(\kappa(r))}\boldsymbol{\varepsilon}(\boldsymbol{u}), \boldsymbol{\varepsilon}(\boldsymbol{v}) )_{\mathcal{H}}
   +(\boldsymbol{\theta}, \boldsymbol{\varepsilon}(\boldsymbol{v}))_{\mathcal{H}} \\
  & \quad \forall\, y=(r, \boldsymbol{\theta})\in Y, \ \forall\, \boldsymbol{u}, \boldsymbol{v} \in V.
  \end{aligned}
  \label{err-10}
\end{equation}
\begin{equation}
\begin{aligned}
(j_c(\gamma_j\boldsymbol{u}), \gamma_j\boldsymbol{v} )_{L^2(\Gamma_3)} = \alpha_j(\gamma_j\boldsymbol{u},\gamma_j\boldsymbol{v})_{L^2(\Gamma_3)} \quad \forall\, \boldsymbol{u}, \boldsymbol{v} \in V.
\end{aligned}
\label{eq:jc specify}
\end{equation}
\begin{equation}
\begin{aligned}
j(\gamma_j\boldsymbol{v}) = \int_{\Gamma_3}j_{\nu}(v_{\nu})d\Gamma \quad \forall\, \boldsymbol{v} \in V.
\end{aligned}
\label{err-13}
\end{equation}
\begin{equation}
  \begin{aligned}
  (\cS u)(t)=\left(\int_0^t\|\boldsymbol{\varepsilon}(\boldsymbol{u}(s))\|_{\mathcal{H}}ds,
  \int_{0}^{t}\mathscr{B}(t-s)\boldsymbol{\varepsilon}(\boldsymbol{u}(s))ds \right)
  \quad \forall\, \boldsymbol{u} \in C(\mathbb{R}_+;V).
  \end{aligned}
  \label{err-11}
\end{equation}
Note that for $j_c$ defined in \eqref{eq:jc specify}, the constants $\alpha_{c}$ and $\alpha_{j}$
in \eqref{pre-jc} are equal: $\alpha_{c}=\alpha_{j}$.

The unique solvability of Problem \ref{Problem 12} has been verified in \cite{SM16}.  Here we consider
fully discrete methods for solving Problem \ref{Problem 12}.   Assume the domain $\Omega$ is
polygonal/polyhedral with a regular family of partitions $\{\mathcal{T}^h\}$.  The linear element space is constructed as follows:
\[V^h=\{ \boldsymbol{v}^h\in C(\overline{\Omega})^d \mid \boldsymbol{v}^h|_{\mathcal{T}} \in \mathcal{P}_1(\mathcal{T})^d
{\rm  \ for \ } \mathcal{T}\in \mathcal{T}^h, \boldsymbol{v}^h=\boldsymbol{0} {\rm \ on \ } \Gamma_1  \}, \]
with $\mathcal{P}_1$ being the space of polynomials of degree no greater than one.
Define
$U^h=\{ \boldsymbol{v}^h\in V^h \mid v_{\nu}^h \leq g {\rm \ at \ node \ points \ on \ } \Gamma_3 \}. $
Assume $g$ is concave; then $U^h \subset U$. Thus the approximation is internal and
the numerical methods for Problem \ref{Problem 12} are defined as follows:

\begin{Problem} \label{Problem contact-expli}
	Find a discrete displacement $\textbf{u}^{kh}:= \{\textbf{u}_n^{kh} \}_{n=0}^N \subset U^h $ such that
		\begin{equation}
		\begin{aligned}
		& (\mathscr{A}\boldsymbol{\varepsilon}(\textbf{u}_n^{kh}), \boldsymbol{\varepsilon}(\textbf{v}^h)-\boldsymbol{\varepsilon}(\textbf{u}_n^{kh}))_{\mathcal{H}}
		+\mu(\boldsymbol{\varepsilon}(\textbf{u}_n^{kh}), \boldsymbol{\varepsilon}(\textbf{v}^h)-\boldsymbol{\varepsilon}(\textbf{u}_n^{kh}))_{\mathcal{H}}\\
		& \qquad\quad -\mu(P_{M(\kappa(\widetilde{\zeta}(t_{n-1})))}\boldsymbol{\varepsilon}(\textbf{u}_{n-1}^{kh}), \boldsymbol{\varepsilon}(\textbf{v}^h)-\boldsymbol{\varepsilon}(\textbf{u}_n^{kh}) )_{\mathcal{H}} \\
		& \qquad\quad +\bigg( \frac{k}{2}\mathscr{B}(t_n-t_0)\boldsymbol{\varepsilon}(\textbf{u}_0^{kh})+k\sum_{j=1}^{n-1}\mathscr{B}(t_n-t_j)\boldsymbol{\varepsilon}(\textbf{u}_j^{kh})\\
		&\qquad\qquad\quad+\frac{k}{2}\mathscr{B}(t_n-t_{n-1})\boldsymbol{\varepsilon}(\textbf{u}_{n-1}^{kh}), \boldsymbol{\varepsilon}(\textbf{v}^h)-\boldsymbol{\varepsilon}(\textbf{u}_n^{kh}) \bigg)_{\mathcal{H}} \\
		& \qquad\quad +\int_{\Gamma_3} j_v^0(u_{n,\nu}^{kh}; v_{\nu}^h-u_{n,\nu}^{kh})d\Gamma+\alpha_j(u_{n,\nu}^{kh}, v_{\nu}^h-u_{n,\nu}^{kh})_{L^2(\Gamma_3)} \\
		&\quad \geq \alpha_j(u_{n-1,\nu}^{kh}, v_{\nu}^h-u_{n,\nu}^{kh})_{L^2(\Gamma_3)} +\langle \boldsymbol{f}_n, \boldsymbol{v}^h-\textbf{u}_n^{kh} \rangle_{V^*\times V}
		\quad \forall\, \textbf{v}^h \in U^h,   \\
		\end{aligned}
		\label{contact-expli}
		\end{equation}
where
\[\widetilde{\zeta}(t_{n-1})=\frac{k}{2}\|\boldsymbol{\varepsilon}(\textbf{u}_0^{kh})\|_{\mathcal{H}}+k\sum_{j=1}^{n-1}
\|\boldsymbol{\varepsilon}(\textbf{u}_j^{kh})\|_{\mathcal{H}}+\frac{k}{2}\|\boldsymbol{\varepsilon}(\textbf{u}_{n-1}^{kh})\|_{\mathcal{H}}.\]	
\end{Problem}

\begin{Problem} \label{Problem contact-impli}
	Find a discrete displacement $\textbf{u}^{kh}:= \{\textbf{u}_n^{kh} \}_{n=0}^N \subset U^h $ such that
		\begin{equation}
		\begin{aligned}
		& (\mathscr{A}\boldsymbol{\varepsilon}(\textbf{u}_n^{kh}), \boldsymbol{\varepsilon}(\textbf{v}^h)-\boldsymbol{\varepsilon}(\textbf{u}_n^{kh}))_{\mathcal{H}}
		+\mu(\boldsymbol{\varepsilon}(\textbf{u}_n^{kh}), \boldsymbol{\varepsilon}(\textbf{v}^h)-\boldsymbol{\varepsilon}(\textbf{u}_n^{kh}))_{\mathcal{H}}\\
		& \qquad\quad -\mu(P_{M(\kappa(\widetilde{\zeta}(t_{n-1})))}\boldsymbol{\varepsilon}(\textbf{u}_{n}^{kh}), \boldsymbol{\varepsilon}(\textbf{v}^h)-\boldsymbol{\varepsilon}(\textbf{u}_n^{kh}) )_{\mathcal{H}} \\
		& \qquad\quad +\bigg( \frac{k}{2}\mathscr{B}(t_n-t_0)\boldsymbol{\varepsilon}(\textbf{u}_0^{kh})+k\sum_{j=1}^{n-1}\mathscr{B}(t_n-t_j)\boldsymbol{\varepsilon}(\textbf{u}_j^{kh})\\
		&\qquad\qquad\quad+\frac{k}{2}\mathscr{B}(t_n-t_{n-1})\boldsymbol{\varepsilon}(\textbf{u}_{n-1}^{kh}), \boldsymbol{\varepsilon}(\textbf{v}^h)-\boldsymbol{\varepsilon}(\textbf{u}_n^{kh}) \bigg)_{\mathcal{H}} \\
		& \qquad\quad +\int_{\Gamma_3} j_v^0(u_{n,\nu}^{kh}; v_{\nu}^h-u_{n,\nu}^{kh})d\Gamma\ge\langle \boldsymbol{f}_n, \boldsymbol{v}^h-\textbf{u}_n^{kh} \rangle_{V^*\times V}
		\quad \forall\, \textbf{v}^h \in U^h.  \\
		\end{aligned}
		\label{contact-impli}
		\end{equation}
\end{Problem}

\begin{Problem} \label{Problem 14}
Find a discrete displacement $\textbf{u}^{kh}:= \{\textbf{u}_n^{kh} \}_{n=0}^N \subset U^h $ such that
\begin{equation}
  \begin{aligned}
   & (\mathscr{A}\boldsymbol{\varepsilon}(\textbf{u}_n^{kh}), \boldsymbol{\varepsilon}(\textbf{v}^h)-\boldsymbol{\varepsilon}(\textbf{u}_n^{kh}))_{\mathcal{H}}
   +\mu(\boldsymbol{\varepsilon}(\textbf{u}_n^{kh}), \boldsymbol{\varepsilon}(\textbf{v}^h)-\boldsymbol{\varepsilon}(\textbf{u}_n^{kh}))_{\mathcal{H}}\\
   & \qquad\quad -\mu(P_{M(\kappa(\widetilde{\zeta}(t_{n-1})))}\boldsymbol{\varepsilon}(2\textbf{u}_{n-1}^{kh}-\textbf{u}_{n-2}^{kh}), \boldsymbol{\varepsilon}(\textbf{v}^h)-\boldsymbol{\varepsilon}(\textbf{u}_n^{kh}) )_{\mathcal{H}} \\
   & \qquad\quad +\bigg( \frac{k}{2}\mathscr{B}(t_n-t_0)\boldsymbol{\varepsilon}(\textbf{u}_0^{kh})+k\sum_{j=1}^{n-1}\mathscr{B}(t_n-t_j)\boldsymbol{\varepsilon}(\textbf{u}_j^{kh})\\
   &\qquad\qquad\quad+\frac{k}{2}\mathscr{B}(t_n-t_{n-1})\boldsymbol{\varepsilon}(\textbf{u}_{n-1}^{kh}), \boldsymbol{\varepsilon}(\textbf{v}^h)-\boldsymbol{\varepsilon}(\textbf{u}_n^{kh}) \bigg)_{\mathcal{H}} \\
   & \qquad\quad +\int_{\Gamma_3} j_v^0(u_{n,\nu}^{kh}; v_{\nu}^h-u_{n,\nu}^{kh})d\Gamma+\alpha_j(u_{n,\nu}^{kh}, v_{\nu}^h-u_{n,\nu}^{kh})_{L^2(\Gamma_3)} \\
  &\quad \geq\alpha_j(2u_{n-1,\nu}^{kh}-u_{n-2,\nu}^{kh}, v_{\nu}^h-u_{n,\nu}^{kh})_{L^2(\Gamma_3)} +\langle \boldsymbol{f}_n, \boldsymbol{v}^h-\textbf{u}_n^{kh} \rangle_{V^*\times V}
   \quad \forall\, \textbf{v}^h \in U^h, n\ge 2,  \\
  \end{aligned}
  \label{err-14}
\end{equation}
and for $n=1$,
\begin{equation}
\begin{aligned}
& (\mathscr{A}\boldsymbol{\varepsilon}(\textbf{u}_1^{kh}), \boldsymbol{\varepsilon}(\textbf{v}^h)-\boldsymbol{\varepsilon}(\textbf{u}_1^{kh}))_{\mathcal{H}}
+\mu(\boldsymbol{\varepsilon}(\textbf{u}_1^{kh}), \boldsymbol{\varepsilon}(\textbf{v}^h)-\boldsymbol{\varepsilon}(\textbf{u}_1^{kh}))_{\mathcal{H}}\\
& \qquad\quad -\mu(P_{M(\kappa(\widetilde{\zeta}(t_{1})))}\boldsymbol{\varepsilon}(\textbf{u}_{1}^{kh}), \boldsymbol{\varepsilon}(\textbf{v}^h)-\boldsymbol{\varepsilon}(\textbf{u}_1^{kh}) )_{\mathcal{H}} \\
& \qquad\quad +\bigg( k\mathscr{B}(t_1-t_0)\boldsymbol{\varepsilon}(\textbf{u}_0^{kh}), \boldsymbol{\varepsilon}(\textbf{v}^h)-\boldsymbol{\varepsilon}(\textbf{u}_1^{kh}) \bigg)_{\mathcal{H}} \\
& \qquad\quad +\int_{\Gamma_3} j_v^0(u_{1,\nu}^{kh}; v_{\nu}^h-u_{1,\nu}^{kh})d\Gamma \geq \langle \boldsymbol{f}_1, \boldsymbol{v}^h-\textbf{u}_1^{kh} \rangle_{V^*\times V}
\quad \forall\, \textbf{v}^h \in U^h.  \\
\end{aligned}
\label{err-14-initial-step}
\end{equation}
\end{Problem}
The numerical scheme for $n=0$ is similar to \eqref{err-14-initial-step} except that the approximation
for the history-dependent term is omitted.

Using arguments similar to that found in \cite{XHHCW19a}, we can show that under the following solution regularity: $\boldsymbol{u} \in W_{loc}^{2,\infty}(\mathbb{R}_+; V)$, $ \boldsymbol{\sigma}\in C(\mathbb{R}_+;H^1(\Omega; \mathbb{S}^d))$, $\boldsymbol{u} \in C(\mathbb{R}_+; H^2(\Omega; \mathbb{R}^d))$, and
$u_{\nu} \in C(\mathbb{R}_+; \widetilde{H}^2(\Gamma_3))$, the following optimal order error bounds hold:
\begin{equation}
  \begin{aligned}
    \max\limits_{0\leq n\leq N}\|\boldsymbol{u}_n-\boldsymbol{u}_{n}^{kh} \|_{V} \leq C(h+ k^{\eta}),
  \end{aligned}
  \label{err-23}
\end{equation}
where $\eta=1$ for Problem \ref{Problem contact-expli} and $\eta=2$ for Problem \ref{Problem contact-impli}, \ref{Problem 14}.

\section{Numerical Results}
\label{sec:numer test}

\begin{figure}[htbp]
	\centering
	\includegraphics[width=6cm]{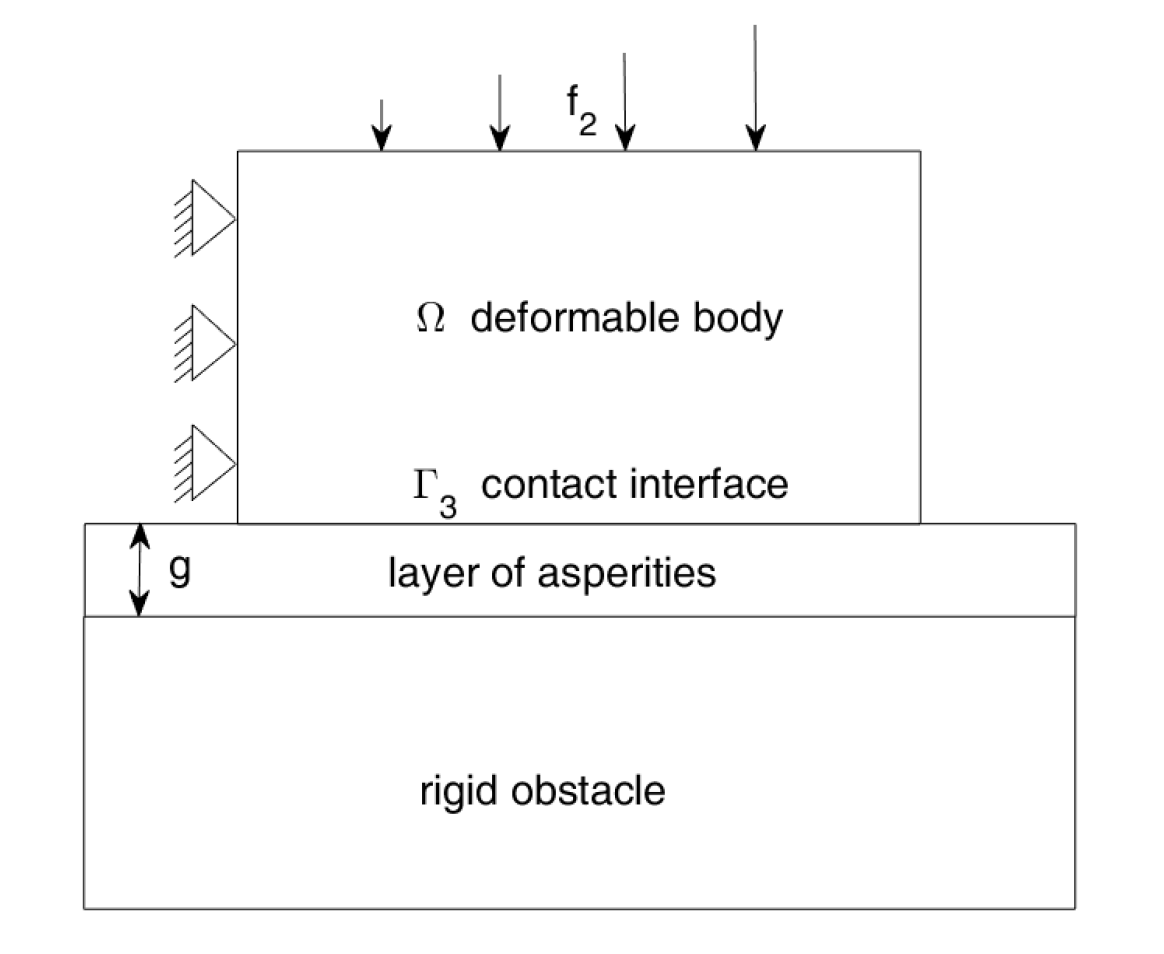}
	\caption{Initial configuration of the contact problem.}
	\label{fig1}
\end{figure}

In this section, we present some numerical results for the three fully discrete schemes stated in
Problem \ref{Problem contact-expli}--\ref{Problem 14}. The same physical setting as depicted
in Figure \ref{fig1} is employed.

Let $\Omega=(0,L_1)\times(0,L_2)$ be a rectangle with boundary $\Gamma$ which is divided into four parts
$$ \Gamma_1=\{0\}\times(0,L_2),\quad \Gamma_2=\{L_1\}\times(0,L_2)\cup[0,L_1]\times \{L_2\},\quad \Gamma_3=[0,L_1]\times\{0\}.$$
For a given $S>0$, the function $j_\nu$ is defined as
\begin{equation}
j_\nu(\xi_\nu) = S\int_{0}^{|\xi_\nu|}\mu_j(s)\,ds
\end{equation}
with
\begin{equation} \label{eq:mu}
\mu_j(s) = \left\{
\begin{aligned}
&0 \qquad s\le0, \\
&c_1s \qquad 0<s\le s_1,\\
&c_1s_1+c_2(s-s_1) \qquad s_1<s\le s_2,\\
&c_1s_2+c_2(s_2-s_1)+c_3(s-s_2)  \qquad s>s_2,
\end{aligned}
\right.
\end{equation}
where $s_1$, $s_2$, $c_1$, $c_2$ and $c_3$ are constants.  The elasticity tensor $\mathscr{A}$ satisfies
\begin{align}
(\mathscr{A}\varepsilon)_{ij} = \frac{E\kappa}{1-\kappa^2}(\varepsilon_{11}+\varepsilon_{22})\delta_{ij}+\frac{E}{1+\kappa}\varepsilon_{ij},
\label{tenorA}
\end{align}
with $1\le i,j\le 2$. $E$ is the Young modulus, $\kappa$ the Poisson ratio of the material and $\delta_{ij}$ denotes the Kronecker symbol.
For the volume and surface forcing, we set
\begin{align}
\boldsymbol{f}_0=(0,-0.1\sin(t) ) \ N/m^2,
\end{align}
\begin{equation}
\boldsymbol{f}_2=\left\{
\begin{aligned}
&(0, 0)  \ N/m   & \text{on} \ \{L_1\} \times (0, L_2), \\
&(0,-0.2\sin(t)\sin(\pi x/2)) \ N/m & \text{on} \ [0,L_1]\times \{L_2\}.
\end{aligned}
\right.
\end{equation}

We test the convergence behavior for the three numerical schemes. The projection on the Von Mises convex is not considered in the convergence tests; thus we let $\mu=0$ in Problem \ref{Problem contact-expli}--\ref{Problem 14}. Values of the other parameters are
\begin{align*}
&L_1=2m,\ L_2=1m,\ E=2N/m^2,\ \kappa=0.3,\\
&\alpha_j=0.5, \ g=0.15m, \ S=1N, \ \mathscr{B}(t) =e^{-t}, \ T=0.5,\\
&s_1 = 0.1, \ s_2 = 0.15, \ c_1 = 0.1, \ c_2 = -0.1, \ c_3 = 0.4.
\end{align*}

The uniform rectangular finite element partitions are introduced to numerically solve the above problem. The numerical solution with $h=k=1/256$ is used as the ``reference" solution in computing numerical solution
errors, and the temporal and spatial convergence orders in the $H^1$ norm will be shown.


\textbf{Example 1 (First order Scheme)} In Tables \ref{1rd_tab1} and \ref{1rd_tab2}, we present the temporal and spatial convergence orders of first-order scheme respectively, and the first-order accuracy in both time and space are shown.

\begin{table}[htbp]
	\centering
	\begin{tabular}{|c|c|c|c|}
		\hline
		$h$ & $k$ & $\|u(\cdot,T)-u_N^{kh}\|_1$ & order  \\
		\hline
		1/256 & 1/4 & 9.82316e-3 & - \\
		\hline
		1/256 & 1/8 &  2.39681e-3 &     2.0351     \\
		\hline
		1/256 & 1/12 & 1.29335e-3 & 1.5215    \\
		\hline
		1/256 & 1/16 & 9.51587e-4 &  1.0667    \\
		\hline
		1/256 & 1/32 & 4.49031e-4 & 1.0835     \\
		\hline
		1/256 & 1/64 &1.93357e-4 &  1.2155 \\
		\hline
	\end{tabular}
	\caption{Convergence orders with spatial step-size fixed for first-order scheme.}
	\label{1rd_tab1}
\end{table}

\begin{table}[htbp]
	\centering
	\begin{tabular}{|c|c|c|c|}
		\hline
		$h$ & $k$ & $\|u(\cdot,T)-u_N^{kh}\|_1$ & order  \\
		\hline
		1/8 & 1/256 & 1.81905e-2 & - \\
		\hline
		1/16 & 1/256 & 1.01388e-2 &  0.8433 \\
		\hline
		1/32 & 1/256 & 5.51935e-3 &   0.8773  \\
		\hline
		1/64 & 1/256 &2.92633e-3 &   0.9154 \\
		\hline
	\end{tabular}
	\caption{Convergence orders with temporal step-size fixed for first-order scheme.}
	\label{1rd_tab2}
\end{table}

\textbf{Example 2 (Second order scheme by fixed-point iteration)} In Tables \ref{2rdFix_tab1} and \ref{2rdFix_tab2}, we present the temporal and spatial convergence orders of second-order fixed-point iteration scheme, respectively, and the second-order accuracy in time, first-order in space are shown.

\begin{figure}[htbp]
	\centering
	\includegraphics[width=10cm]{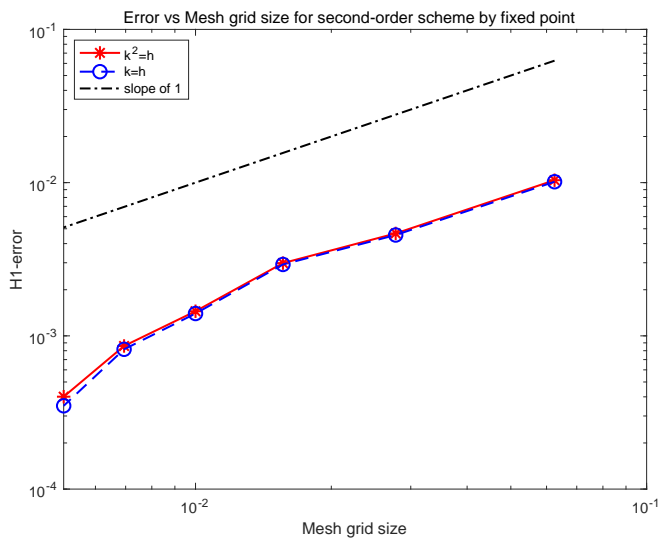}
	\caption{The loglog plot of $H_1$ errors with $h=$1/16, 1/36, 1/64, 1/100, 1/144, 1/196 for second-order fixed-point scheme.}
	\label{fig:2orderFix}
\end{figure}

In addition, we compute the $H_1$ errors for different mesh grid sizes. Two refinement paths are taken to be $k^2=h$ and $k=h$. The results are displayed in Table \ref{Tab03} and the first-order accuracy is shown for both the two refinement paths in Figure \ref{fig:2orderFix}, which indicates the second-order convergence order in time.

\begin{table}[htbp]
	\centering
	\begin{tabular}{|c|c|c|c|}
		\hline
		$h$ & $k$ & $\|u(\cdot,T)-u_N^{kh}\|_1$ & order  \\
		\hline
		1/256 & 1/4 &  2.30136e-3 & - \\
		\hline
		1/256 & 1/8 & 6.06211e-4 &   1.9246    \\
		\hline
		1/256 & 1/12 & 2.75085e-4 &  1.9487    \\
		\hline
		1/256 & 1/16 &1.54881e-4 &   1.9967   \\
		\hline
		1/256 & 1/32 &4.16384e-5 &  1.8952  \\
		\hline
	\end{tabular}
	\caption{Convergence orders with spatial step-size fixed for second-order scheme by fixed-point iteration.}
	\label{2rdFix_tab1}
\end{table}

\begin{table}[htbp]
	\centering
	\begin{tabular}{|c|c|c|c|}
		\hline
		$h$ & $k$ & $\|u(\cdot,T)-u_N^{kh}\|_1$ & order  \\
		\hline
		1/8 & 1/256 & 1.81822e-2 & - \\
		\hline
		1/16 & 1/256 & 1.01334e-2 & 0.8434 \\
		\hline
		1/32 & 1/256 &  5.51570e-3 & 0.8775  \\
		\hline
		1/64 & 1/256 &2.92399e-3 &  0.9156  \\
		\hline
	\end{tabular}
	\caption{Convergence order of the errors with temporal step-size fixed for second-order scheme by fixed-point iteration.}
	\label{2rdFix_tab2}
\end{table}

\begin{table}[htbp]
	\centering
	\begin{tabular}{cccc}
		\toprule
		\multirow{2}{*}{mesh grid size} & \multicolumn{2}{c}{$\|u(\cdot,T)-u_N^{kh}\|_1$} & \multirow{2}{*}{ \shortstack{difference between\\ the front two}}  \\
		\cmidrule(r){2-3}
		&  $k^2=h$      &  $k=h$ \\
		\midrule
		$h=1/16$             & 1.03714e-2  &  1.01335e-2   &  2.38e-4     \\
		$h=1/36$             &  4.66049e-3     & 4.55114e-3  &    1.09e-4   \\
		$h=1/64$             & 2.98351e-3    & 2.92393e-3  &     5.96e-5     \\
		$h=1/100$             & 1.44761e-3    & 1.39882e-3  &     4.88e-5  \\
		$h=1/144$             & 8.58858e-4     & 8.17035e-4  &      4.18e-5    \\
		$h=1/196$             &4.00877e-4      &3.49632e-4   &      5.12e-5     \\
		\bottomrule
	\end{tabular}
	\caption{Comparison of the $H_1$ errors in the refinement path $k^2=h$ and $k=h$ for second-order scheme by fixed-point iteration.}
	\label{Tab03}
\end{table}

\textbf{Example 3 (Second order scheme with extrapolation)} In Tables \ref{2rd_tab1} and \ref{2rd_tab2}, we present the temporal and spatial convergence orders of second-order scheme with extrapolation, respectively, and the second-order accuracy in time, first-order in space are shown.

\begin{figure}[htbp]
	\centering
	\includegraphics[width=10cm]{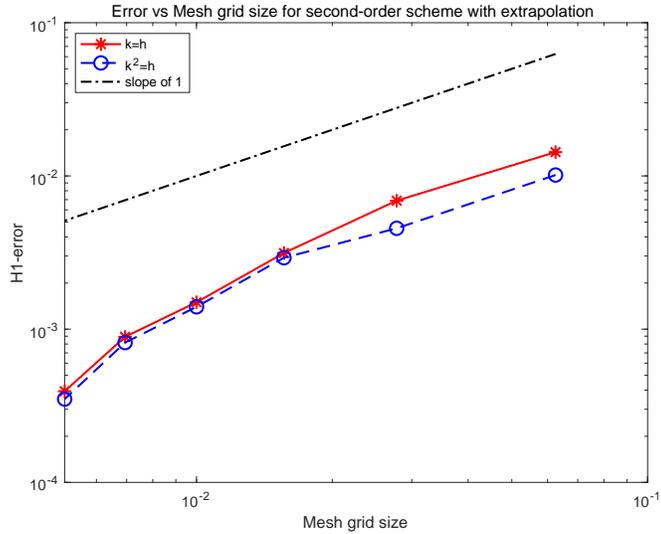}
	\caption{The plot of $H_1$ errors with $h=$1/16, 1/36,  1/64, 1/100, 1/144, 1/196 for second-order scheme with extrapolation.}
	\label{fig:2order}
\end{figure}

In addition, we compute the $H_1$ errors for different mesh grid sizes. Two refinement paths are taken to be $k^2=h$ and $k=h$. The results are displayed in Table \ref{Tab04} and the first-order accuracy is shown for both the two refinement paths in figure \ref{fig:2order}, which indicates the second-order convergence order in time.

\begin{table}[htbp]
	\centering
	\begin{tabular}{|c|c|c|c|}
		\hline
		$h$ & $k$ & $\|u(\cdot,T)-u_N^{kh}\|_1$ & order  \\
		\hline
		1/256 & 1/4 &  1.02222e-2 & - \\
		\hline
		1/256 & 1/8 & 1.15624e-3 &  3.1442  \\
		\hline
		1/256 & 1/12 &3.53015e-4 &    2.9261    \\
		\hline
		1/256 & 1/16 & 2.41930e-4 &   1.3135  \\
		\hline
		1/256 & 1/32 &5.74370e-5 & 2.0745 \\
		\hline
	\end{tabular}
	\caption{Convergence orders with spatial step-size fixed for second-order scheme with extrapolation.}
	\label{2rd_tab1}
\end{table}

\begin{table}[htbp]
	\centering
	\begin{tabular}{|c|c|c|c|}
		\hline
		$h$ & $k$ & $\|u(\cdot,T)-u_N^{kh}\|_1$ & order  \\
		\hline
		1/8 & 1/256 &  1.81823e-2 & - \\
		\hline
		1/16 & 1/256 & 1.01335e-2 & 0.8434    \\
		\hline
		1/32 & 1/256 & 5.51576e-3&  0.8775  \\
		\hline
		1/64 & 1/256 &2.92403e-3 &  0.9156 \\
		\hline
	\end{tabular}
	\caption{Convergence orders with temporal step-size fixed for second-order scheme with extrapolation.}
	\label{2rd_tab2}
\end{table}

\begin{table}[htbp]
	\centering
	\begin{tabular}{cccc}
		\toprule
		\multirow{2}{*}{mesh grid size} & \multicolumn{2}{c}{$\|u(\cdot,T)-u_N^{kh}\|_1$} & \multirow{2}{*}{ \shortstack{difference between\\ the front two}}  \\
		\cmidrule(r){2-3}
		&  $k^2=h$      &  $k=h$ \\
		\midrule
		$h=1/16$             &  1.42988e-2  &   1.01343e-2  &    4.16e-3      \\
		$h=1/36$             &  6.89462e-3  & 4.55108e-3 & 2.34e-3     \\
		$h=1/64$             & 3.13943e-3 & 2.92396e-3  &
		2.15e-4     \\
		$h=1/100$             & 1.49765e-3  & 1.39884e-3   & 9.88e-5     \\
		$h=1/144$             & 8.90577e-4   & 8.17048e-4   & 7.35e-5     \\
		$h=1/196$             &3.93167e-4    &3.49639e-4   &  4.35e-5  \\
		\bottomrule
	\end{tabular}
	\caption{Comparison of the $H_1$ errors in the refinement path $k^2=h$ and $k=h$ for second-order scheme with extrapolation.}
	\label{Tab04}
\end{table}

\begin{figure}[htbp]
	\centering
	\begin{minipage}{4.5cm}
		\includegraphics[width=4.5cm]{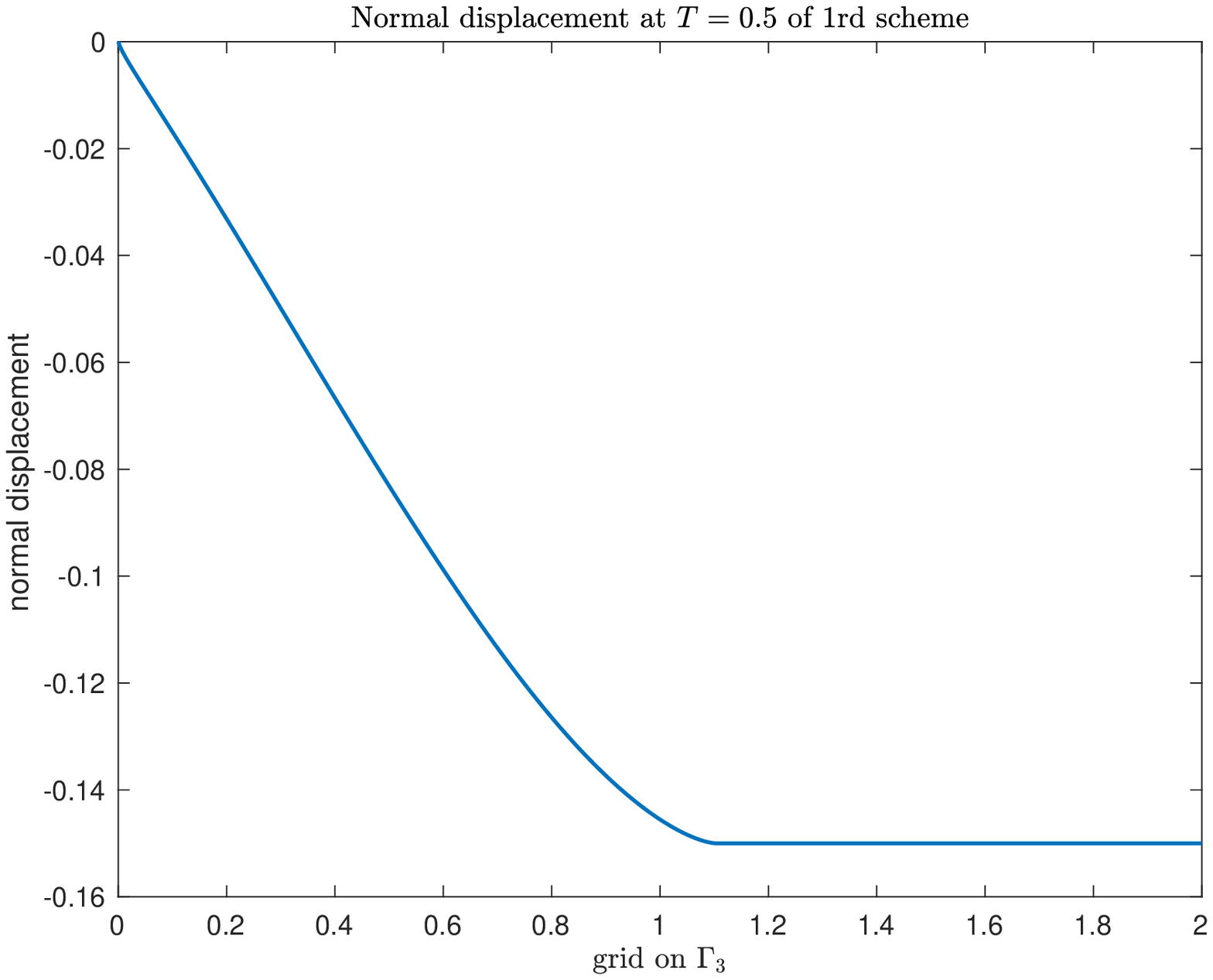}
	\end{minipage}
	\begin{minipage}{4.5cm}
	\includegraphics[width=4.5cm]{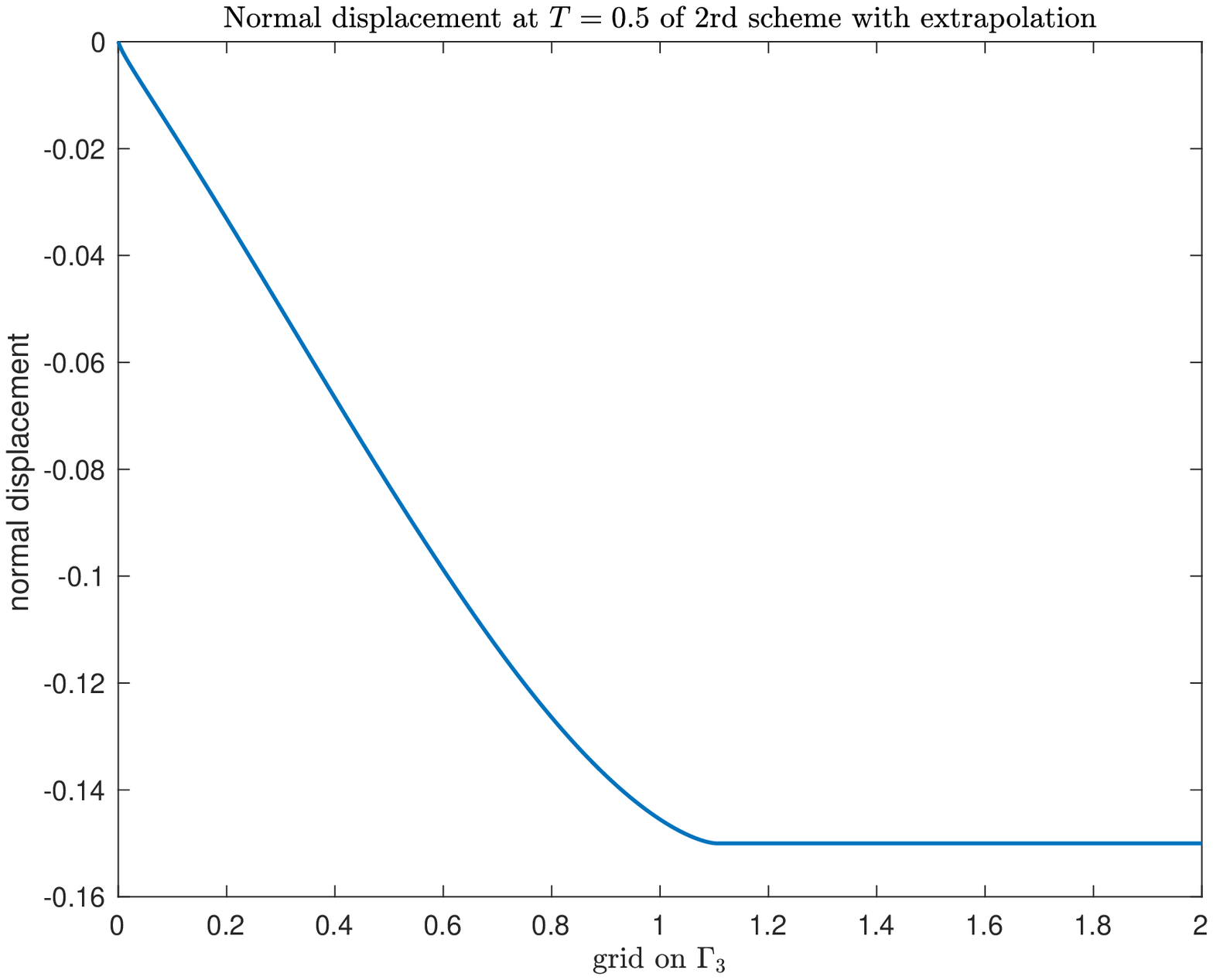}
	\end{minipage}
	\begin{minipage}{4.5cm}
		\includegraphics[width=4.5cm]{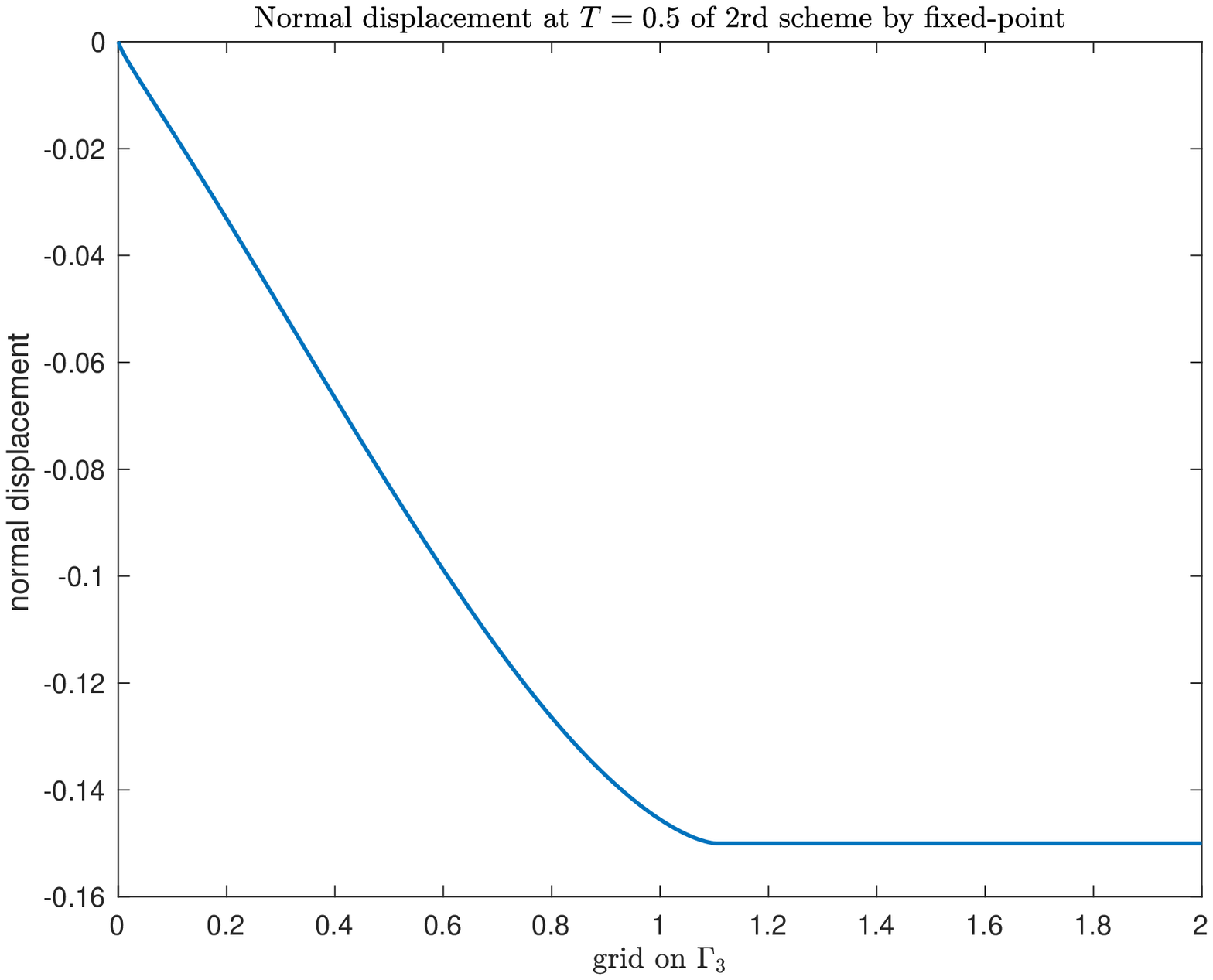}
	\end{minipage}
	\caption{Normal displacement on $\Gamma_3$ at time $T=0.5$ of three numerical schemes.}
	\label{fig:normal displacement with proj}
\end{figure}

In Figure \ref{fig:normal displacement with proj}, the normal displacement on the boundary $\Gamma_3$ at time $T=0.5$ for the three numerical schemes is shown, from which we can see, the maximum penetration is reached as the forcing increased.



\noindent
\textbf{Acknowledgement.}
The work of Wenbin Chen was supported by NSFC under the grants 11671098, 91630309, a
111 Project B08018, and Wenbin Chen also thanks Institute of Scientific Computation and Financial Data Analysis, Shanghai
University of Finance and Economics for the support during his visit.

\end{document}